\documentclass[12pt]{article}
\usepackage[margin=1.1in]{geometry}
\usepackage{subfigure}
\usepackage{amsmath}
\usepackage{amsthm} 
\allowdisplaybreaks
\usepackage{color}
\usepackage{footnote}
\usepackage{algorithm}
\usepackage{algpseudocode}
\usepackage{algorithmicx}
\usepackage{multirow} 
\usepackage{booktabs}
\usepackage{graphicx}
\usepackage{amssymb}
\usepackage{amsbsy}
\usepackage{array}
\usepackage{longtable}
\usepackage{epstopdf}
\usepackage{pbox}
\usepackage{breqn}
\usepackage{mathrsfs}
\usepackage{multicol}
\usepackage{supertabular}
\usepackage{enumerate}
\usepackage[colorinlistoftodos]{todonotes}
\usepackage{bbm}
\usepackage{bbold}
\usepackage{url}
\usepackage[utf8x]{inputenc}
\usepackage[bookmarks]{hyperref}
\usepackage{upgreek}
\usepackage[justification=centering]{caption}

\newtheorem{dfn}{Definition}
\newtheorem{thm}{Theorem}

\newtheorem{prp}{Proposition}

\newtheorem{asm}{Assumption}

\newcommand{\nc}{n^{\mathrm{crit}}}

\newcommand{\Stilde}{\widetilde{F}} 
\newcommand{\vph}{\mathrm{veh/hr}} 
 
\newcommand{\Sbar}{{F}} 
\newcommand{\scrR}{{\mathscr R}} 
\newcommand{\Ftilde}{{\widetilde F}}

\newcommand{\Ntilde}{\widetilde{\mathcal N}}

\newcommand{\ntop}{\overline n} 
\newcommand{\nbot}{\underline n} 
 
\newcommand{\ffspd}{v} 
\newcommand{\conspd}{w} 
\newcommand{\splt}{\beta} 
\newcommand{\inflow}{r} 
\newcommand{\recv}{R} 
\newcommand{\real}{\mathbb{R}} 
\newcommand{\E}{\mathsf{E}} 

\newcommand{\lb}{\left(}
\newcommand{\rb}{\right)}

\newcommand{\I}{\mathcal{I}}

\newcommand{\B}{\mathcal{B}}

\newcommand{\N}{\mathcal{N}}
\newcommand{\R}{\mathcal{R}}

\newcommand{\F}{{\mathcal F}}

\newcommand{\sfp}{\mathsf{p}}

\newcommand{\Ell}{\mathcal{L}}

\newcommand{\scrF}{\mathscr{F}}

\begin{document}
\title{Analysis of a Stochastic Switching Model of Freeway Traffic Incidents}

\author{Li~Jin
        and Saurabh~Amin
\thanks{L. Jin and S. Amin are with the Department of Civil and Environmental Engineering, Massachusetts Institute of Technology, Cambridge, MA 02139 USA (e-mails: jnl@mit.edu; amins@mit.edu)}
}

\maketitle

\begin{abstract}
This article introduces a model for freeway traffic dynamics under stochastic capacity-reducing incidents, and provides insights for freeway incident management by analyzing long-time (stability) properties of the proposed model. Incidents on a multi-cell freeway are modeled by reduction in capacity at the affected freeway sections, which occur and clear according to a Markov chain. We develop conditions under which the traffic queue induced by stochastic incidents is bounded. A necessary condition is that the demand must not exceed the time-average capacity adjusted for spillback. A sufficient condition, in the form of a set of bilinear inequalities, is also established by constructing a Lyapunov function and applying the classical Foster-Lyapunov drift condition. Both conditions can be easily verified for realistic instances of the stochastic incident model. Our analysis relies on the construction of a globally attracting invariant set, and exploits the properties of the traffic flow dynamics. We use our results to analyze the impact of stochastic capacity fluctuation (frequency, intensity, and spatial correlation) on the throughput of a freeway segment. 
\end{abstract}

\section{Introduction}\label{intro}

Freeway traffic networks are prone to capacity disruptions, for example, crashes, road blockage, and other capacity-reducing incidents \cite{jones91,skabardonis97,knoop09,kurzhanskiy10}. In practice, these events can introduce significant congestion in freeways \cite{kwon06,schrank12}. To design traffic control strategies that reduce the congestion and throughput loss resulting from such disruptions, we need to systematically analyze traffic dynamics under stochastic capacity fluctuations. This article introduces a stochastic switching model of freeway traffic dynamics under capacity fluctuations, and studies its stability (in the sense of bounded traffic queue) under fixed inflows.

{\color{black}
According to a recent report by the U.S. Federal Highway Administration \cite[pp. 97]{owens10}, the current practice of managing freeway incidents is to consider each scenario on a case-by-case basis: a traffic manager assesses the impact of detected incidents and responds based on past practice or personal judgment.
Some examples of actually implemented traffic incident management strategies in this class are reported in \cite{carson10}.
This approach is also considered in recent research on traffic management under capacity disruptions \cite{horowitz05,sheu07,long12,kurzhanskiy13}.
}
However, in practice, the implementation of scenario-specific strategies cannot be readily assigned a performance guarantee.
Moreover, the scenario-based approach does not provide guidelines on designing control strategies that are robust to disruptions not captured by the assumed set of scenarios.
Thus, a \emph{model-based approach} for traffic control under disruptions is desirable.
In this article, we follow the model-based approach and study the stability of the well-known cell-transmission model (CTM, see \cite{daganzo94,gomes08}) of freeways under stochastic capacity disruptions.
Our stochastic capacity model captures the main features of capacity-reducing incidents: capacity drop, frequency, and duration.
We combine this model with the CTM and obtain a stochastic switching CTM (i.e. SS-CTM) that can be used to study the throughput and delay in incident-prone freeways.

Before introducing our model and results, we briefly summarize earlier work on stochastic extensions of CTM.
Sumalee et al. \cite{sumalee11} proposed a CTM with parameters subject to random noise; however, their model is largely motivated by the intrinsic randomness in the demand and the flow-density relation, but does not explicitly model capacity-reducing events.
Zhong et al. \cite{zhong14} consider a closely related model in which traffic flow dynamics randomly switch between free-flow and congested modes.
The authors of \cite{zhong14} study an optimal control problem over a finite time horizon.
In contrast, our model admits a richer class of stochastic disruptions.
We focus on the long-time properties of the traffic queue resulting from stochastic disruptions. Indeed, stability of traffic queue is critical for steady-state performance guarantee in any stochastic optimal control of freeway systems, such as that studied in \cite{zhong14}.
Our model is also related to the stochastic models used in \cite{khattak12,miller12} in the context of traffic incident management.
However, these models do not capture the dynamic propagation of incident-induced congestion (spillback), which is a critical mechanism affecting the traffic dynamics in freeways.

In our SS-CTM (formally defined in Section~\ref{Sec_Model}), the capacity of a freeway section switches between a finite set of values (\emph{modes}); the switches are governed by a continuous-time finite-state Markov chain.
The mode transitions can represent abrupt changes in traffic dynamics including:
(i) occurrence of primary incidents \cite{jones91,abdelaty00},
(ii) occurrence of secondary incidents \cite{khattak12,park16}, and
(iii) clearance of incidents \cite{skabardonis97,jin14}.
In a given mode, the traffic density in each section (the continuous state of the model) evolves according to the CTM.
The mode transitions essentially govern the build-up and release of traffic queues in the system.

{\color{black}
We note that actual capacity fluctuation may be more complicated than the finite-state Markov model \cite{jin17,jia10}. However, this model is adequate to study the build-up and release of traffic queues due to stochastic capacity, and also enables a tractable analysis of long-time properties of the traffic queues. 
In a related work \cite{jin17}, we showed that calibration of this model is simple, and that it satisfactorily captures the stochastic variation of capacity in practical situations; see \cite{zhong14,kharoufeh04,baykal09} for related models.
}%

The SS-CTM can be viewed as a \emph{piecewise-deterministic Markov process} (PDMP) without random resets of the continuous state \cite{davis84,benaim15}. 
Since the transition rates are time-invariant, this model is also a specific example of \emph{Markovian jump systems} \cite{boukas07,zhang08}.
In our previous work \cite{jin14}, we provided a preliminary analysis of the accessible set of the model.
In this article, we present a rather complete framework for analyzing the stability of the SS-CTM, based on long-time properties of the deterministic mode dynamics and the Markovian mode transition process.

Our analysis is motivated by the following practical question: on an incident-prone freeway, under what conditions does the long-time average of the total number of vehicles on the freeway (including the upstream queue) remain bounded?
We are interested in deriving intuitive and verifiable stability criteria for the SS-CTM.
Specifically, we consider boundedness of the moment generating function (MGF) of the total number of vehicles on the freeway in our analysis.
In fact, since the traffic density in each cell is upper-bounded by the jam density, boundedness of the total number of vehicles is in line with the boundedness of the upstream queue.
Our notion of stability is similar to that considered by Dai and Meyn \cite{dai95ii}, who considered boundedness of the moments of the queue lengths.

The main results of this article (Theorems~\ref{Thm_Stable1} and \ref{Thm_Stable2}; presented in Section~\ref{Sec_Results}) include a necessary condition and a sufficient condition for stability of the SS-CTM with fixed inflows\footnote{Our setup can be extended to model time-varying or controlled inflows; however, this extension is not within the scope of this article.} and an ergodic mode transition process.
{\color{black}
To establish these results, we build on the theory of stability of continuous-time Markov processes \cite{benaim15,meyn93,gallager13}. Note, however, that the application of standard stability results to our setting is not straightforward due to non-linearity of the CTM dynamics. In this article, we exploit the properties of the mode transition process and the CTM dynamics to develop a tractable approach to characterizing the set of stabilizing inflow vectors -- specifically, an over- and an under-approximation of this set.}
We now introduce the main features of our stability conditions.

Theorem~\ref{Thm_Stable1} is a necessary condition for stability; it states that the SS-CTM is stable only if for each freeway section (cell), the traffic flow does not exceed the time-average of the ``spillback-adjusted'' capacity.
{\color{black}
The notion of spillback-adjusted capacity essentially captures the effect of downstream congestion (i.e. spillback) resulting from capacity fluctuations on the traffic discharging ability of each cell.
To the best of our knowledge, this notion has not been reported previously in the literature. 
The computation of spillback-adjusted capacity builds on the construction of an invariant set that is also globally attracting (Proposition~\ref{Prp_Box}).
Using this construction, we obtain that the capacity fluctuation can result in an unbounded traffic queue even when the inflow in each cell is less than the average capacity; this property essentially results from the effect of traffic spillback.}
Theorem~\ref{Thm_Stable1} also provides a way to identify an over-approximation of the set of stabilizing inflow vectors (i.e. the inflows that satisfy the necessary condition).

{\color{black}
To establish the sufficient condition for stability (Theorem~\ref{Thm_Stable2}), we consider a well-known approach formalized by Meyn and Tweedie \cite{meyn93}.
Application of this approach to our model is challenging, since it requires verification of the Foster-Lyapunov drift condition, i.e. a set of non-linear inequalities, everywhere over the invariant set (which we construct in Proposition~\ref{Prp_Box}).
To resolve this issue, we construct a switched, exponential Lyapunov function that assigns appropriate weights to the traffic densities in various cells, depending on the cells' locations and inflow-capacity ratios.
We also utilize properties of the CTM dynamics to show that the drift condition holds everywhere over the invariant set if a finite set of bilinear inequalities is feasible.
Consequently, standard computational tools \cite{lofberg04} can be used to check whether a given inflow vector maintains stability of SS-CTM or not.}


Overall, these results enable a systematic analysis of the congestion induced by stochastic incidents.
In Section~\ref{Sec_Analysis}, we provide examples to show how Theorems~\ref{Thm_Stable1} and \ref{Thm_Stable2} can be used to characterize the set of stabilizing inflow vectors.
We also comment on the tightness of our results, i.e., the gap between the necessary and the sufficient condition. In addition, we present illustrative examples to study the impact of capacity fluctuation (including its frequency, intensity, and spatial correlation) on the throughput of a freeway section.

Finally, all technical proofs are presented in Appendix.

\section{Stochastic Switching Cell Transmission Model}
\label{Sec_Model}
In this section, we define the stochastic switching cell transmission model (SS-CTM). To develop this model, we introduce a Markovian capacity model, and combine it with the classical CTM \cite{daganzo94}. We also introduce key definitions that are needed for our subsequent analysis.

\subsection{Markovian capacity model}
\label{Sub_Incident}

\begin{figure}[htb]
\centering
\includegraphics[width=0.45\textwidth]{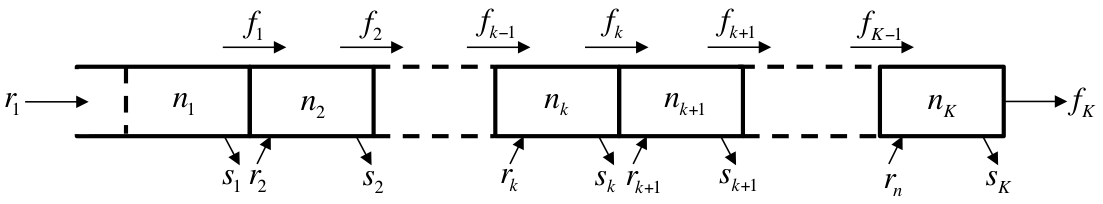}
\caption{SS-CTM with $K$ cells. Cell 1 includes an infinite-sized buffer to accommodate the upstream queue.}
\label{Fig_CTM}
\end{figure}

Consider a freeway consisting of $K$ \emph{cells}, as shown in Figure~\ref{Fig_CTM}. The \emph{capacity} (or saturation rate, in vehicles per hour, or veh/hr) of the $k$-th cell at time $t\ge0$ is denoted by $\Sbar_k(t)$. Let $\Sbar(t)=[\Sbar_1(t),\ldots,\Sbar_K(t)]^T$ denote the vector of cell capacities at time $t$. One can interpret $\Sbar_k(t)$ as the maximum rate at which cell $k$ can discharge traffic to the downstream cell at time $t$.

To model stochastic capacity disruptions, we assume that $\Sbar(t)$ is a finite-state Markov process. Specifically, let $\I$ be a finite set of \emph{modes} of the freeway and let $m=|\I|$. Each mode $i\in\I$ is associated with a vector of cell capacities $\Sbar^i=[\Sbar_1^i,\ldots,\Sbar_K^i]^T$. We define
\begin{subequations}
\begin{align}
\Sbar_k^\min=\min_{i\in\I}\Sbar_{k}^{i},\quad
\Sbar_k^\max=\max_{i\in\I}\Sbar_{k}^{i},\label{Eq_Sbarmax}
\end{align}
\end{subequations}
and refer to $\Sbar_k^\max$ as the \emph{normal} (maximum) \emph{capacity} of cell $k$. For ease of presentation, we assume an identical normal capacity for all cells throughout the article, i.e. $F_k^\max=F^\max$ for $k=1,2,\ldots,K$.

In our model, a mode represents a particular configuration of capacities at various locations (cells). 
We say that the freeway is in the \emph{normal mode} if the maximum capacity is available at every cell. 
We model an incident in cell $k$ by introducing a mode $i$ such that $\Sbar_{k}^i<\Sbar^\max$ and $\Sbar_h^i=\Sbar^\max$ for $h\neq k$. 
Transition from the normal mode to mode $i$ can be viewed as occurrence of an incident in cell $k$; similarly, the transition from $i$ to the normal mode can be viewed as clearance of the incident.
Furthermore, we call the $k$-th cell an \emph{incident hotspot} if $F_k^\min<F^\max$.

Note that the Markovian capacity model can be used to represent more complex situations. For example, two modes can be associated with incidents in the same cell(s), but with different values of capacities, reflecting the difference in incident intensities (e.g. minor and major). 
Furthermore, the occurrence of secondary (or induced) incidents \cite{khattak12} can be modeled as a transition from a mode with an incident in a single cell to a mode with incidents in multiple cells.

{\color{black} Throughout this article, we use $i$ to denote elements of $\I$, and
use $I(t)$ to denote the stochastic mode of the freeway at time $t$.}
The mode $I(t)$ switches according to a continuous-time, finite-state Markov chain defined over the set $\I$ with (time-invariant) transition rates $\{\lambda_{ij};i,j\in\I\}$.
We assume that $\lambda_{ii}=0$ for each $i\in\I$;
note that this is without loss of generality, as inclusion of self-transitions would not affect the traffic flow dynamics.
Let $\nu_i=\sum_{j\in\I}\lambda_{ij}$ and define the \emph{transition matrix} as follows:
\begin{align}
\Lambda=\left[\begin{array}{llll}
-\nu_1 & \lambda_{12} &\ldots & \lambda_{1m}\\
\lambda_{21} & -\nu_2 &\ldots & \lambda_{2m}\\
\vdots & \vdots & \ddots & \vdots\\
\lambda_{m1} & \lambda_{m2} &\ldots & -\nu_m\\
\end{array}\right].
\label{Eq_Lmd}
\end{align}

We assume the following for the mode switching process:

\begin{asm}
\label{Asm_Ergodic}
The finite-state Markov process $\{I(t);t\ge0\}$ is \emph{ergodic}.
\end{asm}

This assumption ensures that the dwell times in each mode are finite almost surely (a.s.). Under this assumption, the process $\{I(t);t\ge0\}$ admits a unique steady-state probability distribution $\sfp=[\sfp_1,\ldots,\sfp_m]$ satisfying:
\begin{align}
\sfp\Lambda=0,\;|\sfp|=1,\;\sfp\ge0,
\label{Eq_pLmd}
\end{align}
where $|\cdot|$ indicates the 1-norm of (row or column) vectors \cite{gallager13}.

\subsection{Traffic flow under stochastic capacities}

The formal definition of the SS-CTM is as follows:

\begin{dfn}
The SS-CTM is a tuple $\langle \I,\N,\R,\Lambda,G\rangle$, where
\begin{itemize}
	\item[-] $\I$ is a finite set of modes (discrete state space) with $|\I|=m$,
	\item[-] $\N=[0,\infty)\times[0,{n}^\max]^{K-1}$ is the set of traffic densities (continuous state space),
	\item[-] $\R\subseteq\real_{\ge0}^{K}$ is the set of inflow vectors,
	\item[-] $\Lambda\in\real^{m\times m}$ is the transition rate matrix governing the mode transitions, and
	\item[-] $G:\I\times\N\times \R\to\real^{K}$ is the vector field governing the continuous dynamics.
\end{itemize}
\end{dfn}

We have defined $\I$ and $\Lambda$, and now introduce $\N$, $\R$, and $G$.

{\color{black}Let ${N}_k(t)$ denote the \emph{traffic density} (in vehicles per mile, veh/mi) in the $k$-th cell at time $t$}, as shown in Figure~\ref{Fig_CTM}. Traffic density ${N}_k(t)$ is non-negative and upper bounded by ${n}^\max_k$, the $k$-th cell's \emph{jam density}.
{\color{black}The $K$-dimensional vector $N(t)=\left[{N}_1(t),{N}_2(t),\dots,{N}_K(t)\right]^T\in\N$ represents the stochastic continuous state of the SS-CTM.}

For ease of presentation, we assume that each cell has the unit length of 1 mi. Furthermore, each cell has identical parameters, including the free-flow speed $v$, the congestion-wave speed $w$, the jam density ${n}^\max$, and -- as mentioned before -- an identical normal capacity $\Sbar^\max$.
The unit of $\ffspd$ and $\conspd$ is miles per hour (mi/hr).
We define the \emph{critical density} of vehicles as
\begin{align}
\nc=\frac{\Sbar^\max}{v}.
\label{Eq_rhoc}
\end{align}
The \emph{sending flows}~$S_k$ and the \emph{receiving flows}~$R_k$ can be written as follows:
\begin{subequations}
\begin{align}
&S_k(i,{n}_k)=\min\left\{\ffspd{n}_k,\Sbar_{k}^{i}\right\},\; k=1,2,\ldots,K,\label{Eq_Sending}\\
&R_k({n}_k)=\conspd({n}^\max-{n}_{k}),\; k=2,3,\ldots,K.\label{Eq_Receiving}
\end{align}
\label{Eq_SR}%
\end{subequations}
Thus, $S_k$ is the traffic flow that cell $k$ can discharge downstream and $R_k$ is the traffic flow from upstream that cell $k$ can accommodate.
The receiving flow of cell 1 will be discussed later in this subsection.
Following \cite{daganzo94}, we assume that,
\begin{align*}
\forall k\in\{1,2,\ldots,K\},\;
\max_{i\in I}S_k(i,\nc)\le R_k(\nc),
\end{align*}
which, along with \eqref{Eq_rhoc} and \eqref{Eq_SR}, implies that,
\begin{align}
\Sbar^\max\le\frac{vw}{v+w}{n}^\max.
\label{Eq_Sbarmax}
\end{align}

Let $r=\left[\inflow_1,\inflow_2,\ldots, \inflow_K\right]^T\in \R=\real_{\ge0}^{K}$ denote the \emph{inflow vector} to the freeway; the unit of $r_k$ is veh/hr. Throughout this article, we assume that the freeway is subject to a fixed (i.e. time-invariant) inflow vector.
Importantly, we also make the standard assumption that, for each cell $k$, the on-ramp flow $r_k$ is prioritized over the sending inflow $S_{k-1}$ from the upstream cell \cite{gomes08}.
Under this priority rule, the flow (in veh/hr) from cell $k$ to cell $k+1$, denoted by $f_k$, is given by the {flow function} (i.e., the so-called~\emph{fundamental diagram}):
\begin{subequations}
\begin{align}
&f_0=0,\\
&f_k(i,{n}_{k},{n}_{k+1},r_{k+1})=
\min\{\beta_kS_k(i,{n}_k),\nonumber\\
&\hspace{0.4in}(\recv_{k+1}({n}_{k+1})- r_{k+1})_+\},\;
k=1,2,\ldots,K-1,\label{Eq_fk}\\
&f_K(i,{n}_K)=\beta_KS_K(i,{n}_K).\label{Eq_fn}
\end{align}
\label{Eq_f(rho)}%
\end{subequations}
where $(\cdot)_+$ stands for the positive part and $\beta_k= f_k/(f_k+s_k)\in(0,1]$ denotes the fixed \emph{mainline ratio}, i.e. the fraction of traffic from cell $k$ entering cell $k+1$. The \emph{off-ramp flow} $s_k$ from cell $k$ is given by
$
s_k(t)=(1/\beta_k-1)f_k(t)
$
for $k=1,2,\ldots,K$.

Let $f(i,{n},r)$ denote the $K$-dimensional vector of flows.
For notational convenience, we denote $S_k(t)=S_k(I(t),{N}_k(t))$, $R_k(t)=R_k({N}_k(t))$, and $f(t)=f(I(t),N(t),r)$.
We say that cell $k$ is experiencing \emph{spillback} at time $t$ if $\beta_kS_k(t)>R_{k+1}(t)-r_{k+1}$, i.e. if the sending flow from cell $k$ exceeds the receiving flow of cell $k+1$.

Due to spillback, there might be traffic queues at the entrances (on-ramps) to the freeway. We track the queue upstream to cell 1 by assuming that cell 1 has a buffer with infinite space to admit this queue, i.e. ${n}^\max_1=\infty$ (see Figure~\ref{Fig_CTM}).
However, we do not consider the on-ramp queues (i.e. queues at on-ramps to cells $2,3,\ldots,K$).
Note that not including the on-ramp queues to cells 2 through $K$ does not affect our stability analysis of the upstream queue, since our priority rule implies that the boundedness of the upstream queue is a sufficient condition for the boundedness of the on-ramp queues.
Hence, we denote the \emph{continuous state space} of the SS-CTM as $\N=[0,\infty)\times[0,{n}^\max]^{K-1}$. 

By mass conservation, traffic density in each cell evolves as follows \cite{gomes08}:
\begin{align}
\dot{N}_k(t)
=f_{k-1}(t)+ r_k-f_k(t)/\splt_k,\; k=1,2,\dots,K.
\label{Eq_rhodot}
\end{align}
From \eqref{Eq_f(rho)} and \eqref{Eq_rhodot}, we can define the vector field $G:\I\times\N\times \R\to\real^n$ governing the continuous state of the SS-CTM as follows:
\begin{subequations}
\begin{align}
&G_1(i,{n},r)=r_1-f_1(i,{n}_1,{n}_2,r_2)/\beta_1,\label{Eq_G1}\\
&G_k(i,{n},r)=f_{k-1}(i,{n}_{k-1},{n}_k,r_k)+r_k\nonumber\\
&\quad-f_k(i,{n}_k,{n}_{k+1},r_{k+1})/\beta_k,\;k=2,3,\ldots,K-1,\label{Eq_Gk}\\
&G_K(i,{n},r)=f_{K-1}(i,{n}_{K-1},{n}_K,r_K)+r_k-f_K(i,{n}_K)/\beta_K.\label{Eq_Gn}
\end{align}
\label{Eq_G}%
\end{subequations}
Note that the vector field $G$ is bounded and continuous in ${n}$.
In a given mode $i\in\I$, the \emph{integral curve} starting from ${n}\in\N$, denoted by $\phi_t^i({n})=[\phi_t^i({n})_1,\ldots,\phi_t^i({n})_K]^T$, can be expressed as follows:
\begin{align}
\phi_t^i({n})={n}+\int_{\tau=0}^tG\Big(i,\phi_\tau^i({n}),r\Big)d\tau.
\label{Eq_phi}
\end{align}

The \emph{hybrid state} of the SS-CTM is $(I(t),N(t))$ at time $t$, and the \emph{hybrid state space} is $\I\times\N$.
The evolution of the discrete (resp. continuous) state is governed by the finite-state Markov process with transition matrix $\Lambda$ (resp.  the vector field $G$).
For an initial condition $(i,{n})\in\I\times\N$, the stochastic process $\{(I(t),N(t));t\ge0\}$ is given by
\begin{subequations}
\begin{align}
&N(t)={n}+\int_{\tau=0}^tG\Big(I(\tau),N(\tau),r\Big)d\tau,\label{Eq_rhot}\\
&\Pr\Big\{I(t+\delta)=j|I(t)=i\Big\}=\lambda_{ij}\delta+\mathrm o(\delta),\; j\neq i.
\end{align}%
\end{subequations}

\subsection{Additional definitions}

For an SS-CTM with a given inflow vector $r\in \R$, the total number of vehicles $|N(t)|$ at time $t$ is given by
\begin{align}
|N(t)|=\sum_{k=1}^K{N}_k(t).
\label{Eq_absrho}
\end{align}
We say that the SS-CTM is \emph{stable} if the moment generating function (MGF) of $N(t)$ is bounded on average; i.e., for some finite constant $C$ and for each initial condition $(i,{n})\in \I\times\N$,
\begin{align}
\limsup_{t\to\infty}\frac1t\int_{\tau=0}^t\E[\exp(|N(\tau)|)]d\tau\le C.
\label{Eq_Stable}
\end{align}
Since a bounded MGF implies a bounded $p$-th moment for all $p\in\mathbb Z_{>0}$, our notion of stability is in line with the notion of bounded moments considered by Dai and Meyn \cite{dai95ii}.
Recall that, in our model, ${N}_2(t),\ldots,{N}_K(t)$ are always upper-bounded by the jam density ${n}^\max$; therefore, $N(t)$ is bounded if and only if ${N}_1(t)$ is bounded.

An alternative notion of stability that is used in the analysis of queueing systems \cite{dai95ii,kulkarni97,jin16} and PDMPs \cite{benaim15,cloez15} is the convergence of the traffic queue towards a unique invariant probability measure.
This notion of stability is equivalent to boundedness of the traffic queue in many simple settings (e.g. for $M/M/1$ queues \cite{gallager13} or fluid queueing systems with stochastic service rates \cite{kulkarni97,jin16}).
However, convergence to a unique invariant probability measure does not always guarantee bounded moments of the traffic queue, which is of practical significance for freeway traffic management.
Therefore, in this article, we consider boundedness of the upstream queue as the stability notion of interest.

A major issue in analyzing the stability of the SS-CTM is to ensure \eqref{Eq_Stable} for all initial conditions $(i,{n})\in \I\times\N$.
We address this issue by constructing a positive invariant set that is also globally attracting and positively invariant \cite{benaim15}, i.e. a set $\Ntilde\subseteq\N$ such that
\begin{subequations}
\begin{align}
&\mbox{(Invariant)}\quad\forall (i,{n})\in \I\times\Ntilde,\; \forall t\ge0,\;
\phi_t^i({n})\in \Ntilde;\label{Eq_Invariant0}\\
&\mbox{(Attracting)}\quad\forall(i,{n})\in \I\times\N,\;
I(0)=i,{N}(0)={n},\nonumber\\
&\hspace{.3in}\forall\epsilon>0,\;\exists T\ge0,\;
\forall t\ge T,\;
\min_{\nu\in\Ntilde}\|N(t)-\nu\|_2\le\epsilon.
\label{Eq_Attracting}
\end{align}
\label{Eq_Invariant}%
\end{subequations}
For convenience, we henceforth refer to any set satisfying \eqref{Eq_Invariant0} and \eqref{Eq_Attracting} simply as an \emph{invariant set}.
Construction of an invariant set considerably simplifies the proofs of our main results (Theorems~\ref{Thm_Stable1} and \ref{Thm_Stable2}), since, to capture long-time behavior of SS-CTM, we only need to consider initial conditions in $\Ntilde$ rather than in $\N$.

Before proceeding further, we introduce two properties of the SS-CTM.
First, the natural filtration $\F_t$ of the SS-CTM is the $\sigma$-algebra generated by $\{(I(s),N(s));s\le t\}$ for all $t\ge0$ \cite{davis84}. Since the realizations of the continuous state are always continuous in time, and since the transition rates $\lambda_{ij}$ are finite and constant, $\F_t$ is \emph{right continuous with left limits} (RCLL, or \emph{c\`adl\`ag} \cite{benaim15}). 
Second, by \cite[Proposition 2.1]{benaim15}, thanks to the RCLL property, the \emph{infinitesimal generator} of the SS-CTM with a fixed inflow $r\in\R$ can be written as an operator $\Ell$ as follows:
\begin{align}
\hspace{-4pt}\Ell g(i,{n})&=G^T(i,{n},r)\nabla_{n} g(i,{n})+\sum_{j\in\I}\lambda_{ij}\Big( g(j,{n})-g(i,{n})\Big),\nonumber\\
&\hspace{1.7in}\forall (i,{n})\in \I\times\N,
\label{Eq_Lg}
\end{align}
where $g:\I\times\N\to\real$ is a function smooth in the second argument, and $\nabla_{n} g(i,{n})$ is the gradient of $g$ with respect to ${n}$.\footnote{We consider $\nabla_{n} g$ as a column vector.}
We utilize the expression of the infinitesimal generator in our stability analysis (while applying the Foster-Lyapunov drift condition \cite{meyn93} in Appendix~C).
\section{Stability of SS-CTM}
\label{Sec_Results}

In this section, we present our results and demonstrate their application via a simple example.
The proofs for these results are available in Appendix.

\subsection{Main results}
\label{Sub_Results}
Our results include a necessary condition (Theorem~\ref{Thm_Stable1}) and a sufficient condition (Theorem~\ref{Thm_Stable2}) for the stability of the SS-CTM under fixed inflows. Both of these conditions rely on the construction of a ``rectangular'' invariant set of the following form:
\begin{align}
\Ntilde=\left[\nbot_1,\infty\right)\times\prod_{k=2}^K\left[\nbot_k,\ntop_k\right].
\label{Eq_Box}
\end{align}

\begin{prp}
\label{Prp_Box}
For an SS-CTM with an inflow vector $r\in\R$, the set $\Ntilde$ of the form in \eqref{Eq_Box} is an invariant set in the sense of \eqref{Eq_Invariant} with the interval boundaries specified as follows:
\begin{subequations}
\begin{align}
&\nbot_1=\min\left\{\frac{r_1}v,\frac{F^\max}{v}\right\},\label{Eq_rhobot1}\\
&\nbot_k=\min\left\{\beta_{k-1}\nbot_{k-1}+\frac{r_{k}}{v},\frac{\beta_{k-1}F_{k-1}^\min+r_k}v,\frac{F^\max}{v}\right\},\nonumber\\
&\hspace{2in} k=2,3,\ldots,K,\label{Eq_rhobotk}\\
&\ntop_K=\left\{\begin{array}{ll}
\dfrac{\beta_{K-1}\Sbar^\max+r_K}v,& \mbox{if }\beta_{K-1}\Sbar^\max+r_K\le\Sbar_K^\min,\\
{n}^\max-\dfrac{\Sbar_{K}^\min}w,&\mbox{o.w.}
\end{array}\right.\label{Eq_rhotopn}\\
&\ntop_k=\left\{\begin{array}{l}
\dfrac{\beta_{k-1}\Sbar^\max+r_k}v,
\quad\mbox{if }\beta_{k-1}\Sbar^\max+r_k\\
\hspace{.47in}\le\min\left\{\Sbar_k^\min,\dfrac{(R_{k+1}(\ntop_{k+1})-r_{k+1})_+}{\beta_{k}}\right\},\\
{n}^\max-\dfrac1w\min\left\{\Sbar_k^\min,\dfrac{(R_{k+1}(\ntop_{k+1})-r_{k+1})_+}{\beta_{k}}\right\},\\
\hspace{2.6in}\mbox{o.w.},
\end{array}\right.\nonumber\\
&\hspace{1.65in} k=K-1,K-2,\ldots,2, \label{Eq_rhotopk}
\end{align}
\label{Eq_rhobottop}%
\end{subequations}
where $S_k$ and $R_k$ are given by \eqref{Eq_SR}.
\end{prp}

The set $\Ntilde$ is constructed by considering the properties of the sending and receiving flows \eqref{Eq_SR}.
Specifically, for each cell $k$, the lower boundary $\nbot_k$ can be viewed as the limiting density when the flow $f_{k-1}$ from upstream is at its minimum and when the flow $f_k$ discharged to downstream is not constrained by cell $(k+1)$'s receiving flow.
Thus, for each $k$ and each ${n}\in\Ntilde$ such that ${n}_k=\nbot_k$, we have $G_k(i,{n},r)\ge0$ in each mode $i\in\I$; i.e. the vector field points in the direction of non-decreasing cell density.
Similarly, on the upper boundary of $\Ntilde$, the vector field points in the direction of non-increasing cell density; i.e. for each $k\ge2$ and each ${n}\in\Ntilde$ such that ${n}_k=\ntop_k$, we have $G_k(i,{n},r)\le0$ in each mode $i\in\I$.

We choose this specific form of invariant set (i.e. Cartesian product of intervals) because of its simple representation \cite{blanchini99}.
Note that, for a given $r\in\R$, Proposition~\ref{Prp_Box} only provides one such construction;
indeed, other rectangular sets satisfying \eqref{Eq_Invariant} exist.
Importantly, this particular construction leads to intuitive and practically relevant conditions (Theorems \ref{Thm_Stable1} and \ref{Thm_Stable2}) that can be used to identify sets of stabilizing and unstabilizing inflow vectors.
In fact, the sharpness of our stability conditions is directly related to the properties of the invariant set; please refer to Proposition~\ref{Prp_Rhotilde} at the end of this subsection.

Before introducing the necessary condition for stability, we need to define a new notion of \emph{spillback-adjusted capacities}: for each $i\in\I$,
\begin{subequations}
\begin{align}
&\Stilde_k^i(\nbot,r)=\min\left\{\Sbar_k^i,\;\frac{1}{\beta_k}\left(R_{k+1}\left(\nbot_{k+1}\right)-r_{k+1}\right)_+\right\},\nonumber\\ 
&\hspace{1.8in}k=1,2,\ldots,K-1,\\
&\Stilde_K^i(\nbot,r)={\Sbar}_K^i,
\end{align}
\label{Eq_Ftilde}%
\end{subequations}
where $\nbot_k$'s are given by \eqref{Eq_rhobot1} and \eqref{Eq_rhobotk}.
Recalling \eqref{Eq_SR} and \eqref{Eq_f(rho)} and noting that $R_k$ is non-increasing in ${n}_k$, one can see that
\begin{align}
&\forall (i,{n})\in \I\times\Ntilde,\;
\forall r\in\R,\;\forall t\ge0,\nonumber\\
&\quad f_k(t)\le\beta_k\min\left\{v{N}_k(t),\Stilde_k^{I(t)}\right\},\;
k=1,2\ldots,K.
\label{Eq_fk<=min}
\end{align}
Thus, $\Stilde_k$ can be interpreted as the capacity adjusted for the receiving flow admissible by the downstream cell $(k+1)$, and hence the name ``spillback-adjusted''.
By considering $\Stilde_k$, we do not need to explicitly involve the receiving flow in our necessary condition for stability.

In addition, we define the following parameters:
\begin{subequations}
\begin{align}
&\beta_k^k=1,\quad k=1,\ldots,K,\\
&\beta_{k_1}^{k_2}=\prod_{h=k_1}^{k_2-1}\beta_h,\quad 1\le k_1\le k_2-1,\; k_2=2,\ldots,K.
\end{align}
\label{Eq_betak1k2}%
\end{subequations}
Note that $\beta_{k_1}^{k_2}$ can be viewed as the fraction of the inflow $r_{k_1}$ that is routed to cell $k_2$.
Thus, for each $k$, we can view $\sum_{h=1}^k\beta_h^kr_h$ as the total \emph{nominal flow} through cell $k$.

Then, we have the following result:

\begin{thm}[Necessary condition]
\label{Thm_Stable1}
Consider an SS-CTM with an inflow vector $r\in\R$.
Let $\sfp_i$ be the solution to \eqref{Eq_pLmd}, $\Stilde_k^i(\nbot,r)$ be as defined in \eqref{Eq_Ftilde}, and $\beta_{k_1}^{k_2}$ be as defined in \eqref{Eq_betak1k2}.
If the SS-CTM is stable in the sense of \eqref{Eq_Stable}, then,
\begin{align}
\sum_{h=1}^k\beta_h^kr_h\le\sum_{i\in\I}\sfp_i\Stilde_k^i\left(\nbot,r\right),\quad
k=1,2,\ldots,K.
\label{Eq_r<=Stilde}
\end{align}
\end{thm}

The left-hand side of \eqref{Eq_r<=Stilde} is the nominal flow through cell $k$. The right-hand side of \eqref{Eq_r<=Stilde} can be viewed as the long-time average of the spillback-adjusted capacity. 
Thus, Theorem~\ref{Thm_Stable1} necessitates that, for the SS-CTM to be stable, the nominal flow cannot exceed the average spillback-adjusted capacity.
This result provides a simple criterion to check for the instability of SS-CTM for a given inflow vector $r\in\R$: if $\sum_{h=1}^k\beta_h^kr_h>\sum_{i\in\I}\sfp_i\Stilde_k^i(\nbot,r)$ for some $k$, then the system is unstable.

An important implication of Theorem~\ref{Thm_Stable1} is that the SS-CTM may be unstable even if, for each cell, the nominal flow is strictly less than the average capacity of the respective cell, i.e.
\begin{align}
\sum_{h=1}^k\beta_h^kr_h<\sum_{i\in\I}\sfp_i\Sbar_k^i,\quad
k=1,2,\ldots,K.
\label{Eq_r<=Stilde2}
\end{align}
To see this, one can note that \eqref{Eq_r<=Stilde2} does not guarantee \eqref{Eq_r<=Stilde}, unless $\Stilde_k^i(\nbot,r)=F_k^i$ for all $k$ and all $i$, which holds only when the inflows are sufficiently small.
In summary, our necessary condition imposes a restriction on the inflow vector that captures the joint effect of capacity fluctuation and spillback.
Note that Theorem~\ref{Thm_Stable1} only involves the steady state probabilities $\sfp_i$ but not the elements of of $\Lambda$ directly.

To develop the sufficient condition, let us limit our attention to the set of inflow vectors satisfying \eqref{Eq_r<=Stilde2}. For each $r$ satisfying \eqref{Eq_r<=Stilde2}, we define the vectors $\gamma=[\gamma_1,\ldots,\gamma_K]^T$ and $\Gamma=[\Gamma_1,\ldots,\Gamma_K]^T$ as follows:
\begin{subequations}
\begin{align}
&\gamma_k=\frac{\sum_{i\in\I}\sfp_i\Sbar_k^i}{\sum_{i\in\I}\sfp_i\Sbar_k^i-\sum_{h=1}^K\beta_h^kr_h},\quad k=1,2,\ldots,K,\\
&\Gamma_K=\gamma_K,\\
&\Gamma_k=\beta_k(\Gamma_{k+1}+\gamma_k),\quad k=K-1,K-2,\ldots,1.
\end{align}
\label{Eq_Wk}%
\end{subequations}
In our sufficient condition, we consider the sum of inflows weighted by $\Gamma_k$:
\begin{align}
\mathscr R(r)=\Gamma^T r.
\label{Eq_scrR}
\end{align}
{\color{black}Essentially, $\Gamma_k$ can be viewed as a weight assigned to the inflow or the traffic density in the $k$-th cell with the following properties: (i) upstream cells have higher weights; (ii) a cell's weight increases with the cell's inflow-capacity ratio.}

In addition, we define the following sets
\begin{subequations}
\begin{align}
&\Theta=\left\{{n}\in\N:{n}_1=\nc,{n}_k\in\{\nbot_k,\ntop_k\},k=2,3,\ldots,K\right\},\label{Eq_Theta}\\
&\hat\Theta=\left\{{n}\in\N:{n}_1=\nbot_1,{n}_k\in\{\nbot_k,\ntop_k\}, k=2,3,\ldots,K\right\},\label{Eq_Thetahat}
\end{align}
\label{Eq_Thetas}%
\end{subequations}
where $\nbot_k$ and $\ntop_k$ are given by \eqref{Eq_rhobottop}.
Note that $\Theta$ and $\hat\Theta$ both have cardinality of $2^{K-1}$, where $K$ is the number of cells.
Furthermore, let
\begin{subequations}
\begin{align}
\mathscr F_i(\nbot,\ntop,r)&=\min_{{n}\in\Theta} \gamma^Tf(i,{n},r),\;i\in\I,\label{Eq_Fi}\\
\hat{\mathscr F}_i(\nbot,\ntop,r)&=\min_{{n}\in\hat\Theta} \gamma^Tf(i,{n},r),\;i\in\I,\label{Eq_Fihat}
\end{align}
\label{Eq_scrF}%
\end{subequations}
where $f(i,{n},r)$ is given by \eqref{Eq_f(rho)}.
Although \eqref{Eq_scrF} involves evaluating minima of $\gamma^Tf$ over discrete sets.
We note that, for typical freeway lengths (in the order of 10 cells), $\mathscr F_i(\nbot,\ntop,r)$ and $\hat{\mathscr F}_i(\nbot,\ntop,r)$ can be obtained by simple enumeration.
As we will show in Appendix~C, $\mathscr F_i(\nbot,\ntop,r)$ and $\hat{\mathscr F}_i(\nbot,\ntop,r)$ can be viewed as lower bounds on the weighted sum of the discharged flows $f_k$ in mode $i$.

Then, we have the following result:

\begin{thm}[Sufficient condition]
\label{Thm_Stable2}
Consider an SS-CTM with an inflow vector $r\in\R$ satisfying \eqref{Eq_r<=Stilde2}. Let $\nbot$ and $\ntop$ be as defined in \eqref{Eq_rhobottop}, $\mathscr R(r)$ as defined in \eqref{Eq_scrR}, and $\mathscr F_i\left(\nbot,\ntop,r\right)$ and $\hat{\mathscr F}_i\left(\nbot,\ntop,r\right)$ as defined in \eqref{Eq_scrF}.

If there exist positive constants $a_{1}$, $a_{2},\ldots$ $a_{m}$ and $b$ such that
\begin{align}
\forall i\in\I,\; a_ib\Big(\mathscr R(r)-\mathscr F_i\left(\nbot,\ntop,r\right)\Big)+\sum_{j\in\I}\lambda_{ij}(a_j-a_i)\le-1,
\label{Eq_BMI}
\end{align}
then, by defining
\begin{subequations}
\begin{align}
&c=\frac1{\max_ia_i},\label{Eq_c}\\
&d\nonumber\\
&=\max_{i\in\I}\Bigg|a_ib\left(\mathscr R(r)-\hat{\mathscr F}_i\left(\nbot,\ntop,r\right)\right)+\sum_{i\in\I}\lambda_{ij}(a_j-a_i)+a_ic\Bigg|\nonumber\\
&\quad\quad\times\exp\Big(b(\Gamma_1{n}_1^{\mathrm{crit}}+\Gamma_2\ntop_2+\cdots+\Gamma_K\ntop_K)\Big),\label{Eq_d}
\end{align}
\label{Eq_cd}
\end{subequations}
we obtain that, for each initial condition $(i,{n})\in\I\times\N$,
\begin{align}
&\limsup_{t\to\infty}\frac1t\int_{\tau=0}^t\E\Big[\exp\Big(N(\tau)\Big)\Big]d\tau\le \left(\frac{d}{c\min_ia_i}\right)^{\frac{1}{b\Gamma_K}}.
\label{Eq_Stable2}
\end{align}
\end{thm}

The bilinear inequalities \eqref{Eq_BMI} essentially restrict the weighted inflow $\scrR$, and thus also restrict the inflow vector $r$ to ensure stability.
The first term on the left-hand side of \eqref{Eq_BMI} captures the difference between the (weighted) inflow and the (weighted) discharged flows;
the second term captures the effect of stochastic mode transitions.
Note that, unlike Theorem~\ref{Thm_Stable1}, Theorem~\ref{Thm_Stable2} explicitly involves the elements of $\Lambda$.

Theorem~\ref{Thm_Stable2} is derived based on an approach introduced by Meyn and Tweedie \cite{meyn93}.
For readers' convenience, we state the relevant result \cite[Theorem 4.3]{meyn93} as follows.
Recall that the SS-CTM is RCLL and its infinitesimal generator $\Ell$ is given by \eqref{Eq_Lg}. Suppose that there exists a norm-like\footnote{The function $V:{\I\times\Ntilde}\to[0,\infty)$ is norm like if $\lim_{{n}\to\infty}V(i,{n})=\infty$ for all $i\in \I$.} function $V:{\I\times\Ntilde}\to\real_{\ge0}$ (called the Lyapunov function) such that, for some $c>0$ and $d<\infty$,
\begin{align}
\Ell V(i,{n})\le -cV(i,{n})+d,\quad\forall (i,{n})\in{\I\times\Ntilde}.
\label{Eq_Drift0}
\end{align}
The above condition is usually referred to as the \emph{drift condition} \cite{meyn93}. Under this condition, for any initial condition $(i,{n})\in \I\times\Ntilde$,
\begin{align}
&\limsup_{t\to\infty}\frac1t\int_{\tau=0}^t\E\Big[V(I(t),{N}(t))|I(0)=i,N(0)={n}\Big]d\tau\nonumber\\
&\hspace{2.5in}\le d/c.
\label{Eq_EV}
\end{align}

We consider the Lyapunov function $V:\I\times\Ntilde\to\real_{\ge0}$ defined as follows:
\begin{align}
V(i,{n})=a_i\exp\lb b\Gamma^T{n} \rb,
\label{Eq_V}
\end{align}
where $a=[a_1,\ldots,a_m]^T$ and $b$ are strictly positive constants (to be determined) and $\Gamma$ is defined in \eqref{Eq_Wk}.
{\color{black}
The switched Lyapunov function captures the effect of both the mode (via the coefficient $a_i$) and the traffic density (via the exponential term $\exp(b\Gamma^T{n})$).
Intuitively, $V$ decreases when the freeway switches to a mode with larger capacities or when traffic is discharged from upstream cells to downstream cells.
 The exponential form is in line with our notion of stability \eqref{Eq_Stable} and facilitates verification of the drift condition \eqref{Eq_Drift0}.}

{\color{black}The main challenge of verifying \eqref{Eq_Drift0} is that 
\begin{align}
d&\ge\max_{(i,{n})\in\I\times\Ntilde}\Ell V(i,{n})+cV(i,{n}).
\label{Eq_LV2}
\end{align}
Since the maximization problem in the right-hand side of \eqref{Eq_LV2} is rather complex, the drift condition is not easy to verify and is thus far from checkable.
To address this challenge, we utilize properties of CTM dynamics to show that \eqref{Eq_Drift0} can be established by minimizing a concave function (see \eqref{Eq_minF} in Appendix C) over the rectangular set $\Ntilde$, where an optimal solution must lie at one of the vertices of $\Ntilde$.
}

{\color{black}
Finally, we note that, in general, there exists a gap between the necessary condition (Theorem~\ref{Thm_Stable1}) and the sufficient condition (Theorem~\ref{Thm_Stable2}).
Since both these results rely on the invariant set $\Ntilde$, the gap depends on the construction of the invariant set.
Indeed, both results also apply to other invariant sets that can be expressed of the form in \eqref{Eq_Box}.
The following result addresses how the construction of the invariant set affects the gap:

\begin{prp}
Consider two invariant sets
\begin{align*}
\Ntilde=[\nbot_1,\infty)\times\prod_{k=2}^K[\nbot_k,\ntop_k],\; \Ntilde'=[\nbot_1',\infty)\times\prod_{k=2}^K[\nbot_k',\ntop_k']
\end{align*}
such that $\Ntilde\subseteq\Ntilde'$, i.e. $\nbot\ge\nbot'$ and $\ntop\le\ntop'$. For a given $r\in\R$,
\begin{enumerate}[(i)]
\item if $\Ftilde_k^i(\nbot,r)$ satisfies \eqref{Eq_r<=Stilde}, so does $\Ftilde_k^i(\nbot',r)$;
\item if $\scrF_i(\nbot',\ntop',r)$ allows positive solutions for $a_1,\ldots,a_m$, $b$ to \eqref{Eq_BMI}, so does $\scrF_i(\nbot,\ntop,r)$.
\end{enumerate}
\label{Prp_Rhotilde}
\end{prp}

Proposition \ref{Prp_Rhotilde} implies that a smaller invariant set leads to sharper stability conditions (i.e. a smaller gap between the necessary condition and the sufficient condition).
Indeed, the invariant set given by Proposition~\ref{Prp_Box} is in some cases the smallest invariant set of the form in \eqref{Eq_Box};
see the next subsection for examples.

In Section~\ref{Sub_r}, we show via an example that the gap typically decreases as the model becomes ``simpler'' (i.e. fewer incident hotspots, fewer modes, smaller extent of spillback, etc.), and that, for some cases, the gap vanishes.
Other factors affecting the gap include estimate of the spillback-adjusted capacity and the form of the Lyapunov function. Addressing these issues is important from a practical viewpoint and is part of our ongoing work.
}

\subsection{Example: stability of a two-cell, two-mode SS-CTM}
\label{Sub_Two1}

\begin{figure}[hbt]
\centering
\includegraphics[width=0.24\textwidth]{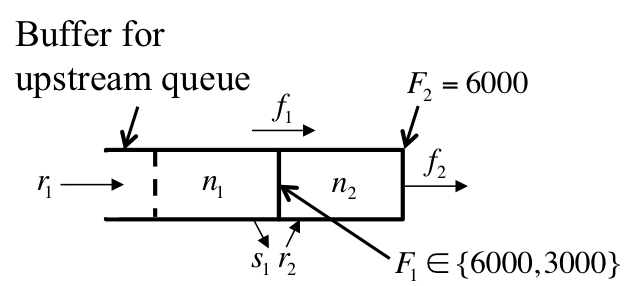}
\caption{A two-cell freeway with a single incident hotspot.}
\label{Fig_Two1}
\end{figure}

Now we use our results to study the stability of a two-cell SS-CTM as in Figure~\ref{Fig_Two1}. The deterministic version of this system has been studied by Kurzhanskiy \cite{kurzhanskiy09}.

We consider the following stochastic capacity model:
\begin{align}
&\I=\{1,2\},
\Sbar^1=[6000,6000]^T,\Sbar^2=[3000,6000]^T,\nonumber\\
&\lambda_{12}=\lambda_{21}=1[hr^{-1}].
\label{Eq_Incident}
\end{align}
In words, the SS-CTM has a nominal capacity of 6000 veh/hr (mode 1). Cell 1 is an incident hotspot, whose capacity reduces to 3000 veh/hr during incidents (mode 2). Incidents occur and clear at the rate 1 per hour. 
The mainline ratios are $\beta_1=0.75$ and $\beta_2=1$. The inflow is $r=[r_1,r_2]^T\in \R=\real_{\ge0}^2$.
Other model parameters are given by Table~\ref{Tab_Parameters}. 
The steady-state distribution of the mode transition process can be easily obtained:
$
\sfp_1=\sfp_2=0.5.
$

\begin{table}[hbt]
\centering
\caption{Traffic flow parameters in example SS-CTM.}
\begin{tabular}{|l|c|c|c|}
\hline
Name                      & Symbol    & Value & unit                \\ \hline
Cell length               & $l$         & 1 & mi                      \\ \hline
Free-flow speed           & $v$         & 60 & mi/hr                    \\ \hline
Congestion wave speed     & $w$         & 20 & mi/hr                    \\ \hline
Jam density               & ${n}^\max$      & 400 & veh/mi               \\ \hline
Critical density          & $\nc$        & 100 & veh/mi            \\ \hline
\end{tabular}
\label{Tab_Parameters}
\end{table}

Our objective is to determine whether the system is stable under given inflow vectors $r=[r_1,r_2]^T\in R$.
Specifically, we consider two different cases, one being unstable and the other being stable.

\emph{1)} Consider the inflow vector $r=[4320,2400]^T$. 

\begin{figure}[hbt]
\centering
\includegraphics[width=0.35\textwidth]{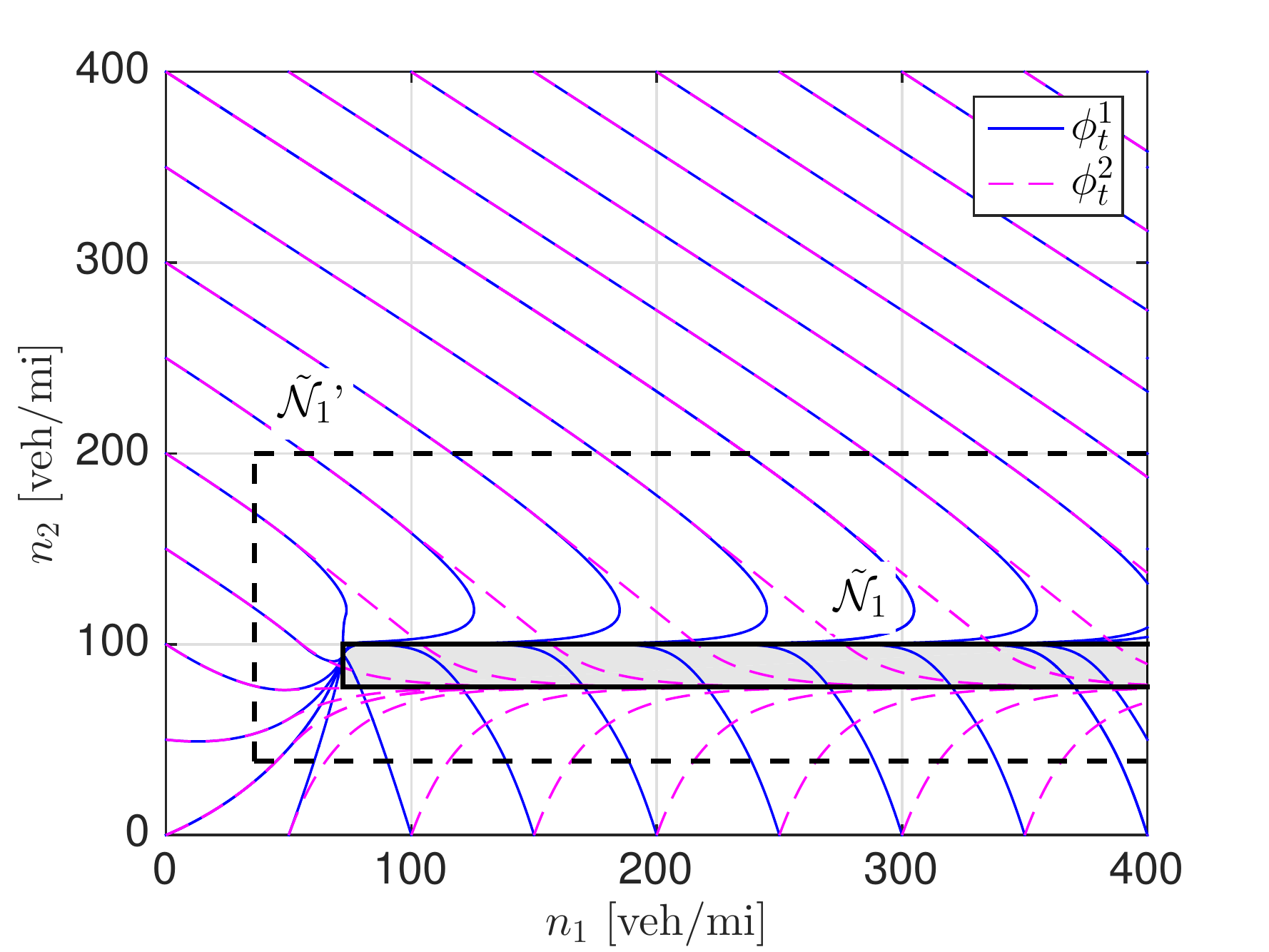}
\caption{Integral curves and invariant set (shaded area) of the two-cell SS-CTM with $r=[4320,2400]^T$.}
\label{Fig_ASet1}
\end{figure}

Note that the nominal flows given by this $r$ do not exceed the average capacities of both cells, i.e. satisfy \eqref{Eq_r<=Stilde}:
\begin{align*}
&r_1=4320<0.5(6000)+0.5(3000)=4500[\vph],\\
&0.75r_1+r_2=0.75(4320)+2400)=5640<6000[\vph].
\end{align*}
However, na\"ively concluding stability based on the above observation is incorrect.
By computing the invariant set and the spillback-adjusted capacities, and applying Theorem~\ref{Thm_Stable1}, we can conclude that the SS-CTM is in fact unstable.
To see this, we first use Proposition~\ref{Prp_Box} to obtain an invariant set
\begin{align}
\Ntilde_1=[72,\infty)\times[77.5,100].
\label{Eq_Rhotilde1}
\end{align}
Figure~\ref{Fig_ASet1} illustrates that, for the two-dimensional system, all integral curves are attracted to $\Ntilde_1$ and remain within it thereafter (i.e. $\Ntilde_1$ satisfies \eqref{Eq_Invariant}).
In fact, it turns out that $\Ntilde_1$ is the smallest invariant set of the form in \eqref{Eq_Box}.
Thus, we have
\begin{align*}
&\Stilde_1^1(\nbot,r)=4960[\vph],\\
&\Stilde_2^1(\nbot,r)=3000[\vph],\\
&\Stilde_2^1(\nbot,r)=\Stilde_2^2(\nbot,r)=6000.
\end{align*}
We can check that
\begin{align*}
r_1=4320>\sfp_1\Stilde_1^1+\sfp_2\Stilde_1^2
=3980.
\end{align*}
Therefore, by Theorem~\ref{Thm_Stable1}, the system is unstable. 

Now, we use this example to show that a smaller invariant set leads to a sharper necessary condition. Consider another set $\N'=[36,\infty)\times[38.75,200]$; see Figure~\ref{Fig_ASet1} to note that this set is also invariant in the sense of \eqref{Eq_Invariant}. 
For $\Ntilde'_1$, we have
\begin{align*}
\Stilde_1^1(\nbot',r)=\Stilde_2^1(\nbot',r)=\Stilde_2^2(\nbot',r)=6000,\;\Ftilde_1^2(\nbot',r)=3000.
\end{align*}
Since $r_1=4320>0.5(6000)+0.5(3000)=4500$, $\beta_1r_1+r_2=5640\le6000$, the given inflow vector satisfies the necessary condition associated with $\Ntilde'_1$.
In other words, using the condition with $\Ntilde'_1$, we cannot conclude the instability of the given inflow vector $r$, while using the condition with $\Ntilde_1$ (as given by \eqref{Eq_Rhotilde1}) we are able to do so.

\emph{2)} Now, consider the inflow vector $r=[3600,600]^T$.

\begin{figure}[hbt]
\centering
\includegraphics[width=0.35\textwidth]{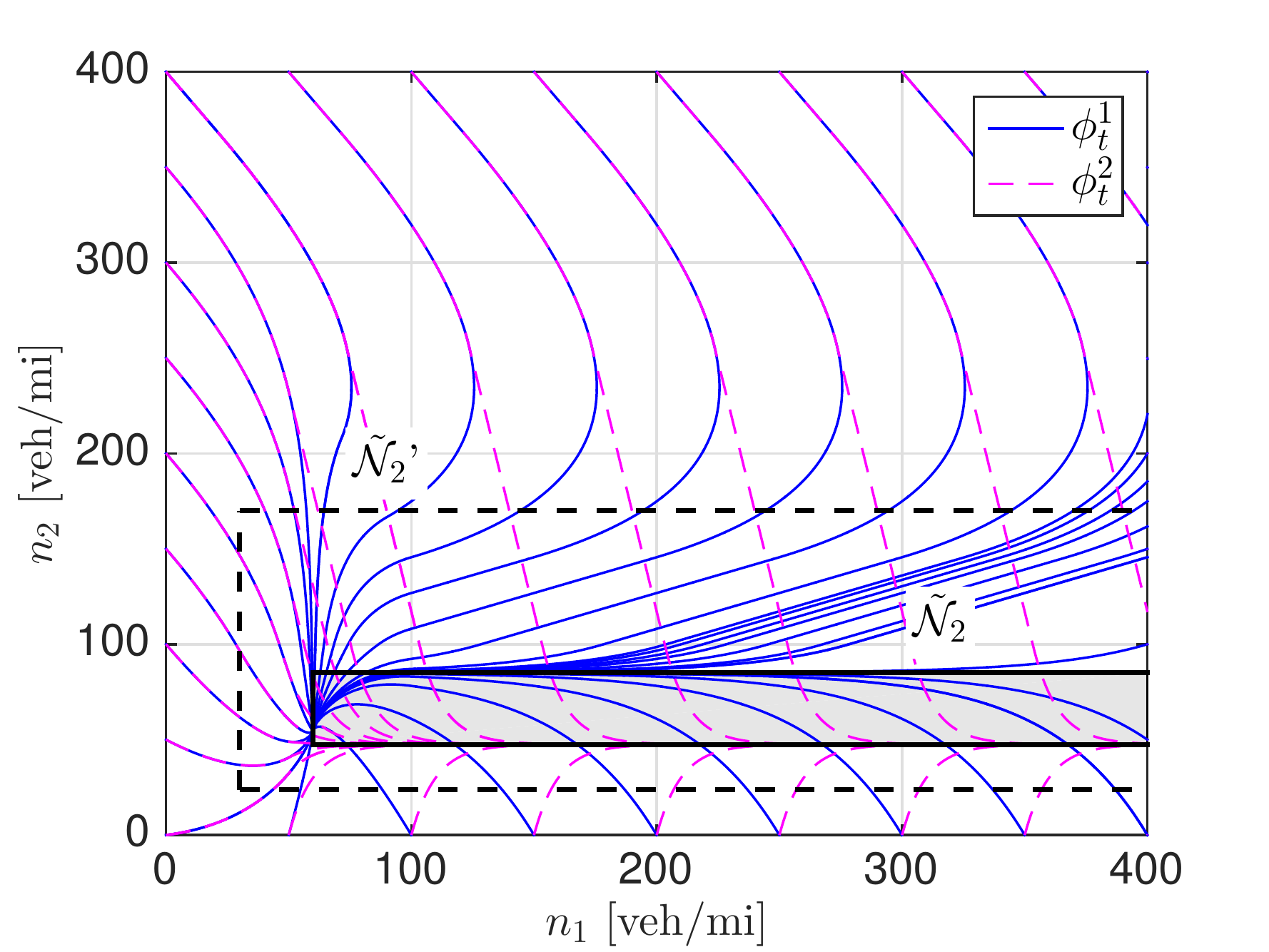}
\caption{Integral curves and invariant set (shaded area) of the two-cell SS-CTM with $r=[3600,600]^T$.}
\label{Fig_ASet2}
\end{figure}

Once again, from Proposition~\ref{Prp_Box}, we obtain the following invariant set:
\begin{align}
\Ntilde_2=[60,\infty)\times[47.5,85].
\label{Eq_Rhotilde2}
\end{align}
Similar to the previous case, this invariant set is the smallest one that can be written as a Cartesian product of intervals; see Figure~\ref{Fig_ASet2}.
Then, we obtain from \eqref{Eq_Ftilde} that
\begin{align*}
&\Stilde_1^1(\nbot,r)=6000[\vph],\;
\Stilde_1^2(\nbot,r)=3000[\vph],\\
&\Stilde_2^1(\nbot,r)=\Stilde_2^2(\nbot,r)=6000[\vph].
\end{align*}
Note that, the necessary condition (Theorem~\ref{Thm_Stable1}) is satisfied, since
\begin{align*}
&r_1\le0.5(6000)+0.5(3000)=4500[\vph],\\
&\beta_1r_1+r_2=3300\le6000[\vph],
\end{align*}

To check the sufficient condition (Theorem~\ref{Thm_Stable2}), we first compute the weights using \eqref{Eq_Wk}:
\begin{align*}
\gamma_1=5,\;\gamma_2=2.22,\;
\Gamma_1=5.42,\;\Gamma_2=2.22,
\end{align*}
which implies
\begin{align*}
&\mathscr R(r)=2.08\times10^4,\\
&\scrF_1(\nbot,\ntop,r)=2.88\times10^4,\;
\scrF_2(\nbot,\ntop,r)=1.76\times10^4,
\end{align*}
Then, plugging the numerical values into \eqref{Eq_BMI}, we obtain
\begin{align*}
&a_1b(2.08-2.88)\times10^4+(a_2-a_1)\le-1,\\
&a_2b(2.08-1.76)\times10^4+(a_1-a_2)\le-1.
\end{align*}
We can easily obtain a solution to the above bilinear inequalities, one of which is
\begin{align*}
a_1=10,\;a_2=17,\;b=10^{-4}.
\end{align*}
Hence, we conclude from Theorem~\ref{Thm_Stable2} that the system is stable in the sense of \eqref{Eq_Stable}.

Now we show that an invariant set larger than $\Ntilde_2$ (as given by \eqref{Eq_Rhotilde2}) leads to a weaker sufficient condition.
Consider a larger set $\Ntilde'_2=[35,\infty)\times[23.75,170]$ (see Figure~\ref{Fig_ASet2}), which is also invariant in the sense of \eqref{Eq_Invariant}.
With $\Ntilde'_2$, we have
\begin{align*}
&\scrR(r)=2.08\times10^{4},\\
&\scrF_1(\nbot',\ntop',r)=2.57\times10^4,\;
\scrF_2 (\nbot',\ntop',r)=1.44\times10^4,
\end{align*}
which leads to the sufficient condition that
\begin{align*}
&a_1b(-0.49)10^4+(a_2-a_1)\le-1,\\
&a_2b(0.64)10^4+(a_1-a_2)\le-1.
\end{align*}
It is not hard to show that the above inequalities have no positive feasible solutions for $a_1$, $a_2$, and $b$. Thus, the sufficient condition with $\Ntilde'_2$ cannot be used to ensure stability.
Hence, the sufficient condition resulting from using $\Ntilde_2$ is sharper than that obtained by using $\Ntilde'_2$.
\section{Some Practical Insights}
\label{Sec_Analysis}

In this section, we derive some practical insights for freeway traffic management under capacity fluctuations.
Specifically, we use our results to
(i) identify the set of stable inflow vectors for a given capacity model,
(ii) analyze the impact due to frequency and intensity of capacity fluctuation on throughput, and
(iii) study the effect of correlation between capacity fluctuation at various locations.

\begin{figure}[hbt]
\centering
\includegraphics[width=0.24\textwidth]{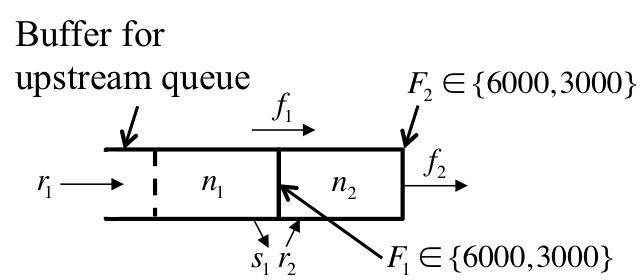}
\caption{A two-cell freeway with two incident hotspots.}
\label{Fig_Two2}
\end{figure}

\begin{figure*}
\centering
\subfigure[Baseline model.]{
\centering
\includegraphics[width=0.3\textwidth]{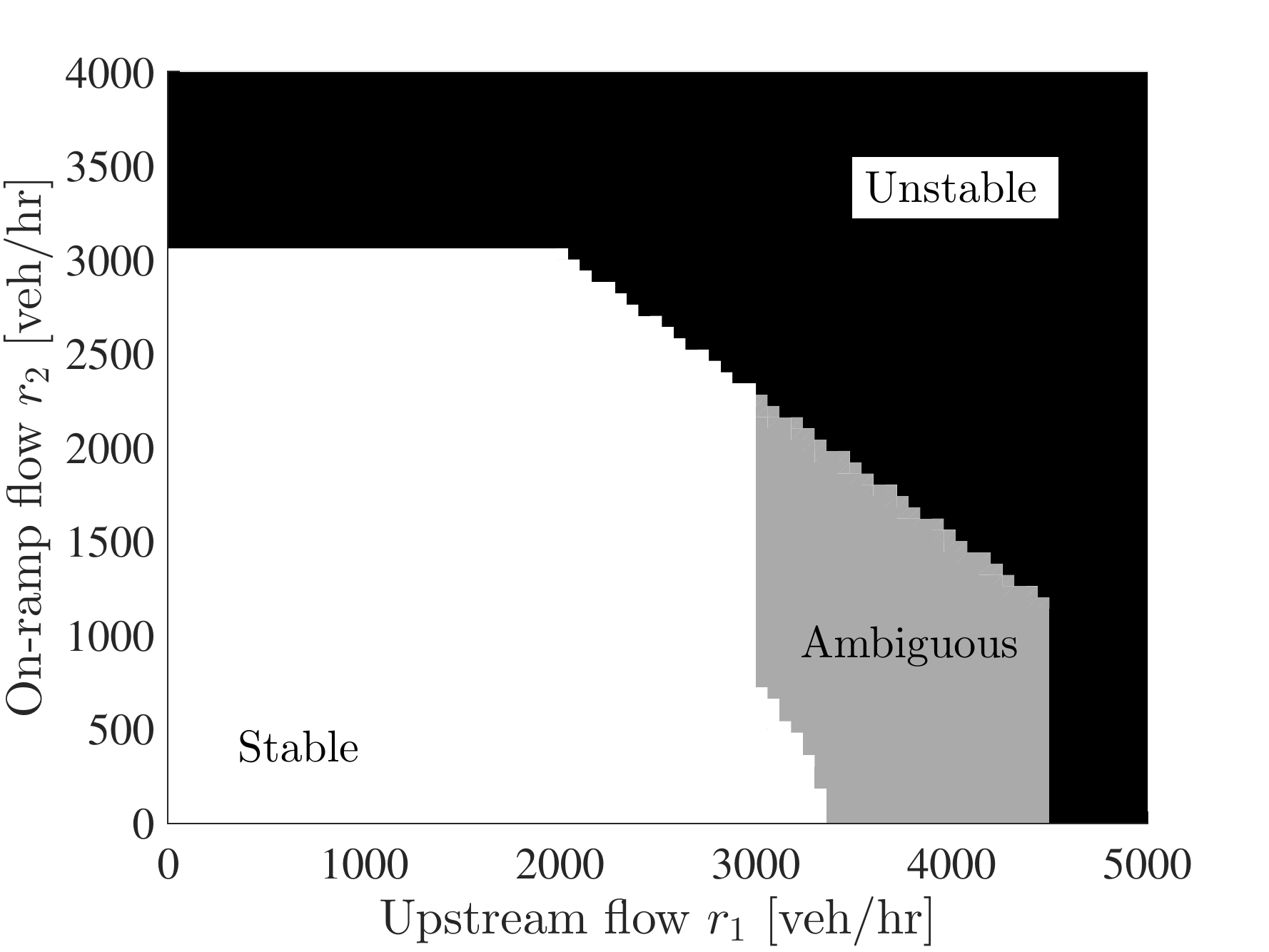}
\label{Fig_Stable0}
}
\subfigure[Variant 1.]{
\centering
\includegraphics[width=0.3\textwidth]{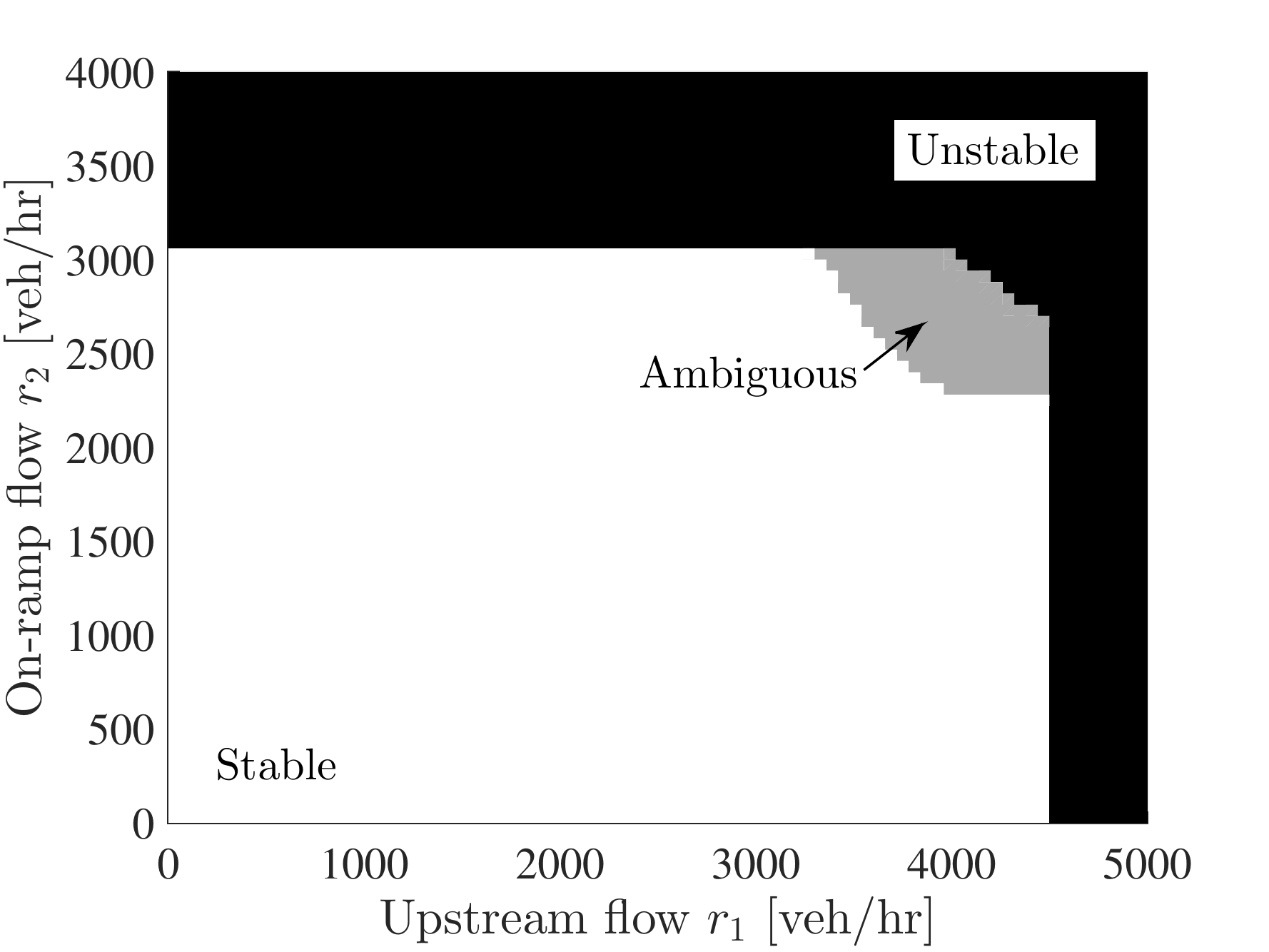}
\label{fig_gap}
}
\subfigure[Variant 2.]{
\centering
\includegraphics[width=0.3\textwidth]{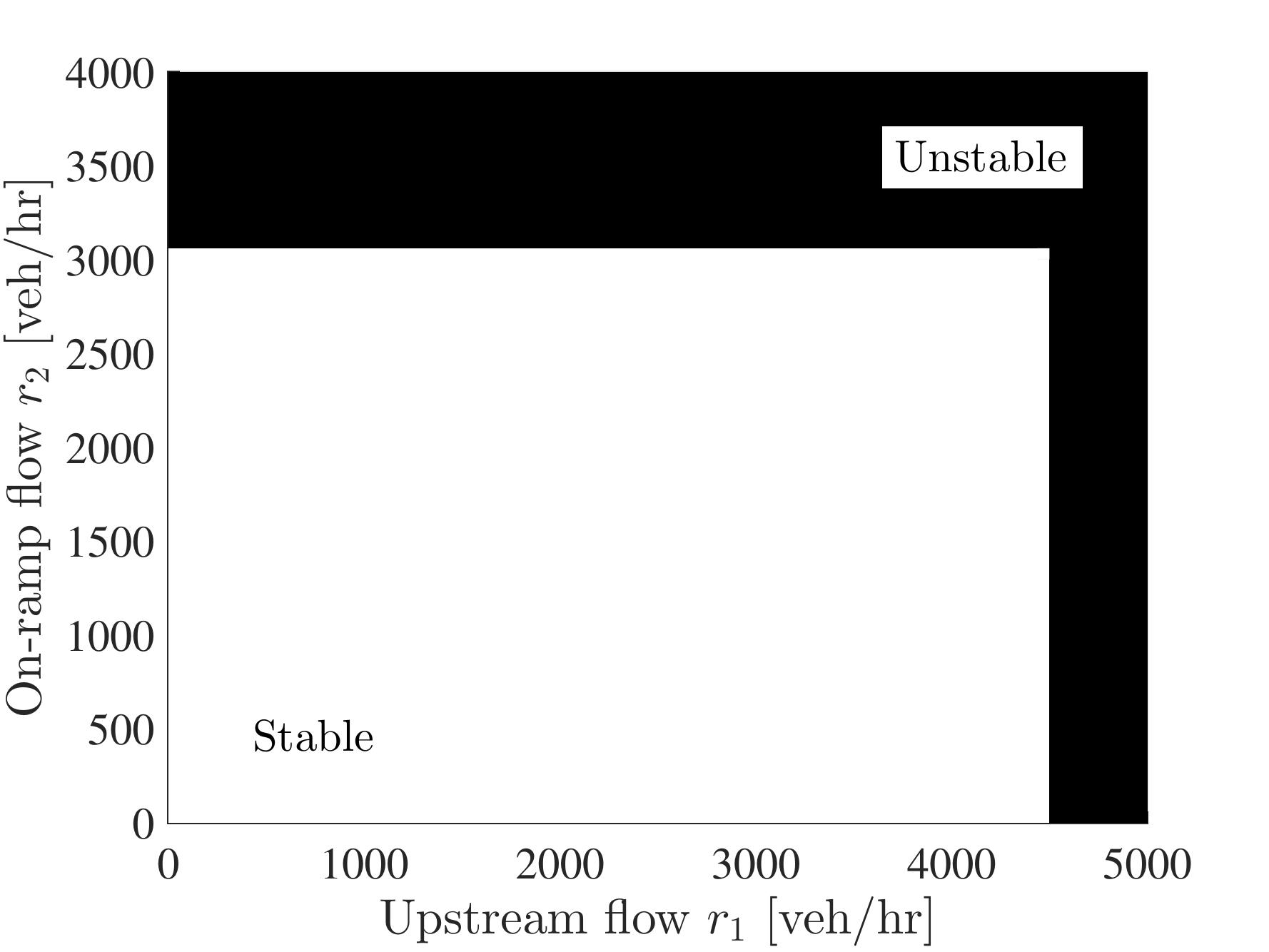}
\label{fig_nogap}
}
\caption{Stability of the two-cell freeway with various inflow vectors $r=[r_1,r_2]^T$ determined by Theorems \ref{Thm_Stable1} and \ref{Thm_Stable2}.}
\label{Fig_gap}
\end{figure*}

A two-cell system as shown in Figure~\ref{Fig_Two2} is sufficient for our purpose.
Note that this system is distinct from that considered in Section~\ref{Sub_Two1} in that both cells of this system have stochastic capacities. 
The parameters for the (baseline) capacity model is as follows:
\begin{subequations}
\begin{align}
&\I=\{1,2,3,4\}.\\
&\Sbar^1=[6000,6000]^T,\;
\Sbar^2=[3000,6000]^T,\label{eq_F12}\\
&\Sbar^3=[6000,3000]^T,\;
\Sbar^4=[3000,3000]^T,\label{eq_F34}\\
&\Lambda=\left[\begin{array}{cccc}
-2&1&1&0\\
1&-2&0&1\\
1&0&-2&1\\
0&1&1&-2
\end{array}\right].
\end{align}
\label{eq_model}%
\end{subequations}
Note that the transition rate matrix defined above implies that the capacity fluctuations at both cells are mutually independent.
(We will discuss about the impact of correlation between cell capacities later.)
The remaining parameters for the SS-CTM are given by Table~\ref{Tab_Parameters}.

\subsection{Set of stabilizing inflow vectors}
\label{Sub_r}

For an inflow vector $r=[r_1,r_2]^T\in\real_{\ge0}^2$, we know that $r$ is unstable if it does not satisfy the necessary condition (Theorem~\ref{Thm_Stable1}), and that $r$ is stable if it satisfies the sufficient condition (Theorem~\ref{Thm_Stable2}).
For practicality, we also assume that the on-ramp has a fixed saturation rate of 3000 veh/hr. Thus, the on-ramp inflow $r_2$ cannot exceed 3000 veh/hr if the freeway is stable.

Applying our stability conditions to various inflow vectors $[r_1,r_2]^T$, we obtain Figure~\ref{Fig_Stable0}.
In this figure, the $r_1$-$r_2$ plane is partitioned into three regions:
The ``Unstable'' region (in black) depicts the set of inflow vectors violating the necessary condition.
Thus, any inflow vector in this region leads to an infinite traffic queue.
We denote the complement of this region as $\R_1$.
The ``Stable'' region (in white) depicts the set of inflow vectors satisfying sufficient condition.
In this example, we solve the bilinear inequalities \eqref{Eq_BMI} using YALMIP, a MATLAB-based optimization tool \cite{lofberg04}. 
By Theorem~\ref{Thm_Stable2}, the inflow vectors in this region lead to a traffic queue bounded in the sense of \eqref{Eq_Stable}.
We denote this region as $\R_2$.

Notice that there is a gap, labeled as ``Ambiguous'', between the ``Stable'' and the ``Unstable'' regions.
This region shows the gap between the necessary condition and the sufficient condition; for inflow vectors in this region, our stability conditions do not provide a conclusive answer.


These results can be used to calculate stabilizing inflows that lead to maximum throughput.
For a stable inflow vector $r$, we follow \cite{kurzhanskiy10} and define throughput $J(r)$ of the two-cell system as follows:
\begin{align*}
J(r)=r_1(2l)+r_2l=2lr_1+lr_2.
\end{align*}
One can interpret $J(r)$ as the average \emph{vehicle miles traveled} per unit time.
Let $\R_{\mathrm{stable}}$ be the set of stabilizing inflow vectors.
For the given capacity model, we are interested in the maximum throughput $J_\max$:
\begin{align*}
J_\max=\max_{r\in \R_{\mathrm{stable}}}2lr_1+lr_2
\end{align*}
Since our stability conditions are in general not necessary and sufficient, it is not easy to compute the exact value of $J_\max$. However, we are able to derive upper and lower bounds for $J_\max$.
Note that $\R_2\subseteq \R_{\mathrm{stable}}\subseteq \R_1$.
Hence, we can obtain an upper bound $\overline J_\max$ and a lower bound $\underline J_\max$ for $J_\max$ as follows:
\begin{align*}
&\overline J_\max=\max_{r\in \R_1}(2lr_1+lr_2),\\
&\underline J_\max=\max_{r\in \R_2}(2lr_1+lr_2).
\end{align*}
For this particular example, we can numerically obtain the following bounds:
\begin{align}
7170\le J_\max\le 8910[\mbox{veh-mi/hr}].
\label{Eq_Jmax1}
\end{align}

Next, to illustrate how model characteristics affect the size of the ``Ambiguous'' region, we consider the following variants of the baseline model specified by \eqref{eq_model}:
\begin{enumerate}
\item Cell 2 is no longer an incident hotspot; i.e. $F_2^i=6000$ veh/hr for all $i\in\{1,2,3,4\}$.
\item Cell 2 is no longer an incident hotspot, and its capacity is doubled; i.e. $F_2^i=12000$ veh/hr for all $i\in\{1,2,3,4\}$.
\end{enumerate}
Compared to the baseline model, both variants has fewer incident hotspots; furthermore, variant 2 has a sufficiently large capacity of cell 2 to prevent spillback between cells 1 and 2.
Hence, we can sort the three models in the order of ``decreasing complexity'' as follows: baseline, variant 1, variant 2.

By comparing Figures~\ref{Fig_Stable0}--\ref{fig_nogap}, we can see that the gap between the sufficient condition and the necessary condition reduces as the model complexity reduces.
In addition, for variant 2, as shown in Figure~\ref{fig_nogap}, the ambiguous region vanishes; i.e. the two-cell freeway is stable if and only if the inflows are less than average capacities ($r_1<4500$ veh/hr, $r_2<3000$ veh/hr.)\footnote{In fact, for this simple example, one can prove this statement by explicitly constructing $a_i$ and $b$ in the Lyapunov function under this condition.}

\subsection{Impact of capacity fluctuation on throughput}
\label{Sub_DeltaSbar}

Now we estimate the maximum throughput that can be achieved under a class of capacity models parameterized by $\Delta F$ and $\lambda$ as follows:
\begin{align*}
&\I=\{1,2,3,4\}.\\
&\Sbar^1=[6000,6000]^T,\;
\Sbar^2=[6000-\Delta\Sbar,6000]^T,\\
&\Sbar^3=[6000,6000-\Delta\Sbar]^T,\\
&\Sbar^4=[6000-\Delta\Sbar,6000-\Delta\Sbar]^T,\\
&\Lambda=\left[\begin{array}{cccc}
-2\lambda&\lambda&\lambda&0\\
1&-(1+\lambda)&0&\lambda\\
1&0&-(1+\lambda)&\lambda\\
0&1&1&-2
\end{array}\right].
\end{align*}
By varying the parameters $\Delta F$ and $\lambda$, we can study the effect the intensity and the frequency of capacity fluctuations on the maximum throughput of SS-CTM.

\begin{figure}[hbt]
\centering
\subfigure[Throughput decreases as incident frequency increases, with $\Delta F$ fixed at $3000$ veh/hr.]{
\centering
\includegraphics[width=0.35\textwidth]{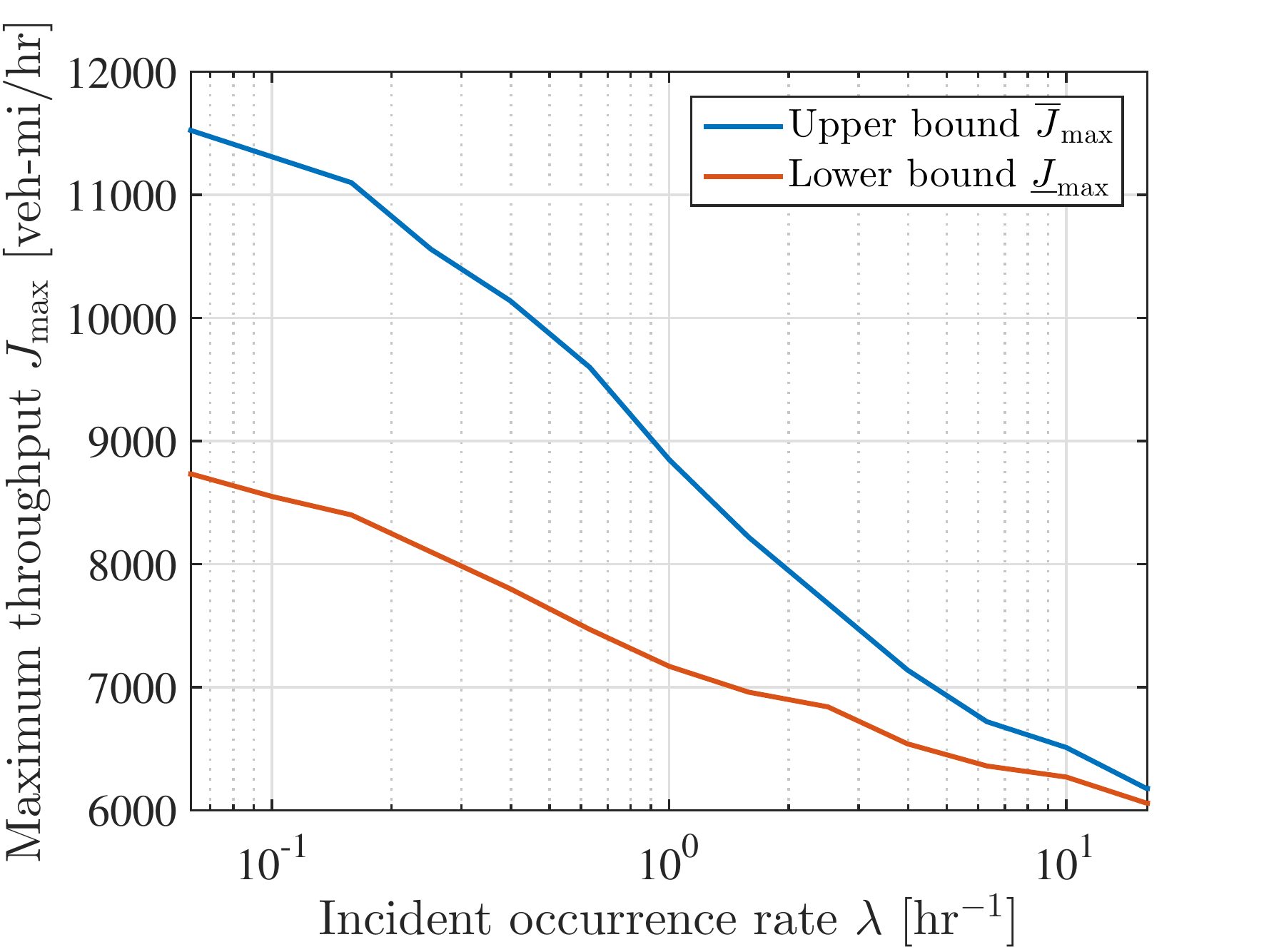}
\label{Fig_Jvlambda}
}
\subfigure[Throughput decreases as incident intensity increases, with $\lambda$ fixed at 1 per hour.]{
\centering
\includegraphics[width=0.35\textwidth]{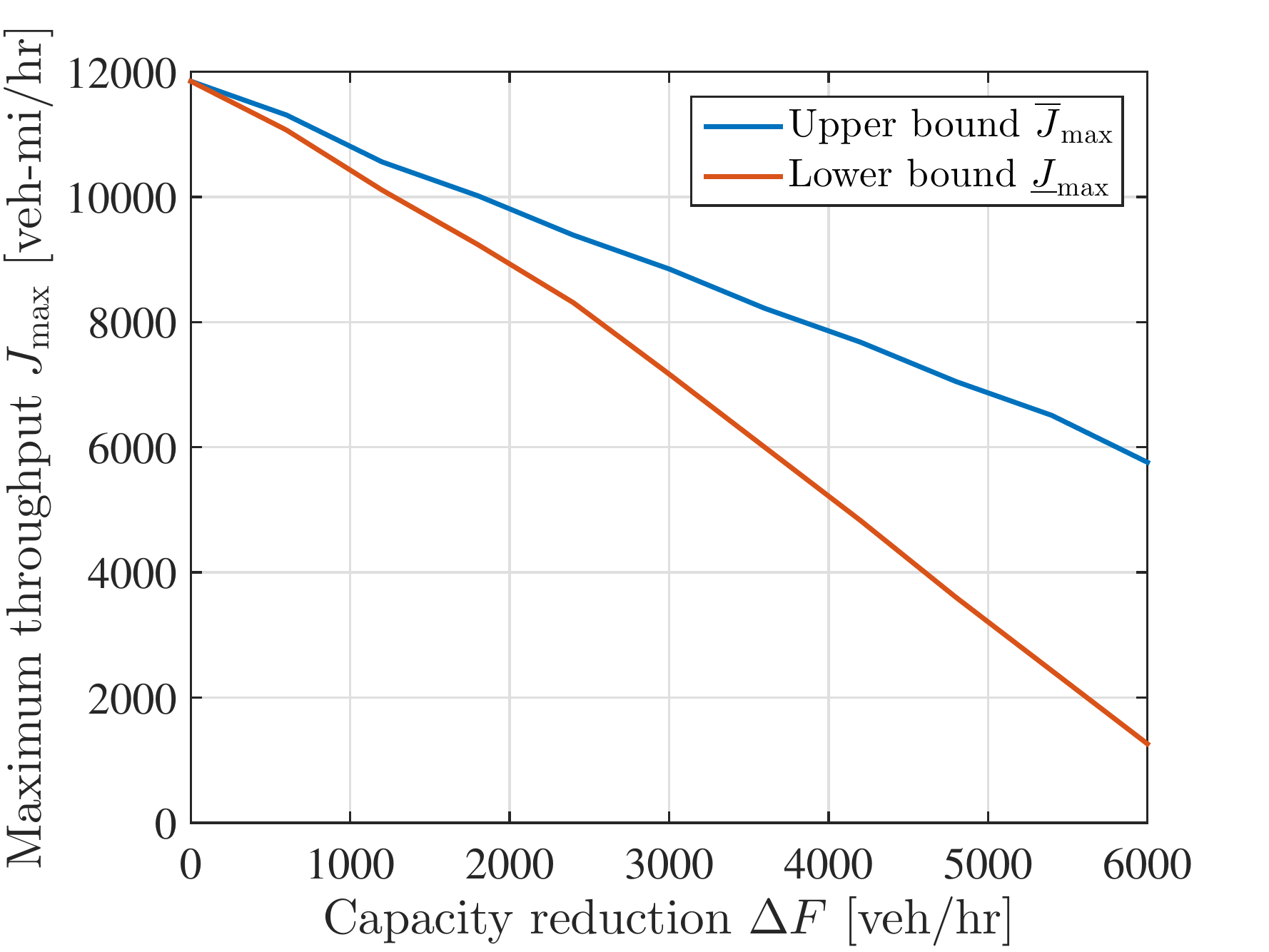}
\label{Fig_JvDeltaSbar}
}
\caption{Relation between maximum throughput and incident frequency/intensity. The upper (resp. lower) bounds result from Theorem~\ref{Thm_Stable1} (resp. Theorem~\ref{Thm_Stable2}).}
\label{Fig_J}
\end{figure}

For various $(\lambda,\Delta\Sbar)$ pairs, we numerically determine $\overline J_\max$ and $\underline J_\max$. Figure~\ref{Fig_Jvlambda} shows that, with $\Delta\Sbar$ fixed at 3000 veh/hr, $\overline J_\max$ and $\underline J_\max$ decreases as $\lambda$ increases. This is intuitive: more frequent capacity disruptions leads to lower throughput.
Figure~\ref{Fig_JvDeltaSbar} shows that, with $\lambda$ fixed at 1 per hr, $\overline J_\max$ and $\underline J_\max$ decreases as $\Delta\Sbar$ increases. This is also intuitive: larger capacity reduction leads to lower throughput.

Note that the throughput tends to be more sensitive to $\Delta\Sbar$ than to $\lambda$. 
Indeed, as shown in Figure~\ref{Fig_Jvlambda}, if $\lambda$ is doubled from 1 (the baseline) to 2 per hr, the upper (resp. lower) bound is reduced by 7\% (resp. 3\%). 
However, if $\Delta\Sbar$ is doubled from 3000 (the baseline) to 6000 veh/hr, we can observe from Figure~\ref{Fig_JvDeltaSbar} that the upper (resp. lower) bound is reduced by 35\% (resp. 82\%). 

The gap between $\overline J_\max$ and $\underline J_\max$ in Figures~\ref{Fig_Jvlambda} and \ref{Fig_JvDeltaSbar} result from the gap between the necessary condition and the sufficient condition for stability.

\subsection{Impact of correlated capacity fluctuation}
\label{Sub_Lambda}

So far we have assumed that the cell capacities in the two-cell freeway are independent. Now we consider the case where the cells' capacities are correlated.
We consider two extreme cases as follows.

\emph{Case 1}: Suppose that the capacities of both cells are identical for all $t\ge0$. In other words, the freeway has two modes $\{1,2\}$ associated with $\Sbar^1=[6000,6000]^T$ and $\Sbar^2=[3000,3000]^T$. In this case, the transition matrix is
\begin{align*}
\Lambda_{\mathrm{Case\,1}}=\left[\begin{array}{cc}-1&1\\1&-1\end{array}\right].
\end{align*}
By implementing Theorems~\ref{Thm_Stable1} and \ref{Thm_Stable2}, we obtain the following bounds for the maximum throughput:
\begin{align*}
7485\le J_\max\le8910\mbox{[veh-mi/hr]}.
\end{align*}
The upper bounds is the same as that in \eqref{Eq_Jmax1}, since the necessary condition \eqref{Eq_r<=Stilde} only involves the spillback-adjusted capacities of individual cells but does not depend on the correlation between cells. However, the lower bound (obtained from the sufficient condition) is higher compared to that in \eqref{Eq_Jmax1}.

\emph{Case 2}: Suppose that the freeway always has exactly one incident in either cell. In other words, the freeway has two modes $\{1,2\}$ associated with $\Sbar^1=[6000,3000]^T$ and $\Sbar^2=[6000,3000]^T$. In this case, the transition matrix is
\begin{align*}
\Lambda_{\mathrm{Case\,2}}=\left[\begin{array}{cc}-1&1\\1&-1\end{array}\right].
\end{align*}
By implementing Theorems~\ref{Thm_Stable1} and \ref{Thm_Stable2}, we obtain the following bounds for the maximum throughput:
\begin{align*}
6720\le J_\max\le8910\mbox{[veh-mi/hr]}.
\end{align*}
In this case, the lower bound is lower than that in the benchmark case \eqref{Eq_Jmax1}.

In conclusion, even with the same cell average cell capacities (4500 veh/hr for both cells), the freeway's maximum throughput estimated by the sufficient condition could vary (in the order of 10\%) due to spatial correlation.
Therefore, compared to traditional approaches that assume independent cell capacities (e.g. \cite{sumalee11}), our approach is able to capture the effect of spatial correlation and possibly achieve better throughput.
In a related work \cite{jin17}, we have found that the extent of spatial correlation can be quite significant in practice.

{We finally note that, in some situations, the transition rate matrix $\Lambda$ may not be easy to calibrate.
Analysis of the SS-CTM with partially known transition rates is indeed a practically relevant question, but is beyond the scope of this paper.
We refer the readers to \cite{zhang08} for more information on stochastic switching systems with partially known transition probabilities.}
\section{Concluding Remarks}
\label{Sec_Concluding}

The main contributions of this article are: (i) a stochastic switching cell transmission model for traffic dynamics in incident-prone freeways and (ii) easily checkable stability conditions for the SS-CTM (Theorems~\ref{Thm_Stable1} and \ref{Thm_Stable2}). 
A sufficient condition for stability is that the inflow does not exceed the spillback-adjusted capacity.
A necessary condition for stability is that a set of bilinear inequalities, which is derived from the Foster-Lyapunov drift condition, admit positive solutions.
Both conditions build on a construction of a globally attracting and invariant set of the SS-CTM (Proposition~\ref{Prp_Box}).
Using these results, we derive new insights for the impact of capacity fluctuation on the upstream traffic queue length and the attainable throughput.

Our work can be extended in several directions. First, with a different (or more sophisticated) choice of Lyapunov function or a tighter invariant set, the stability conditions may be further refined. Second, estimation methods for the SS-CTM can be developed based on previous work such as \cite{jin17,dervisoglu09}. Third, recent results on network traffic flow models \cite{jin16,coogan15,como13i} makes possible the extension of our model and method to the network setting, and to feedback-controlled systems.
\appendix
\label{Sec_Stability}

\section*{Appendix\\Proofs of Main Results}
\subsection{Proof of Proposition~\ref{Prp_Box}}




\subsubsection{Invariant} To show the invariance of $\Ntilde$, we demonstrate that the vector field $G$ points towards the interior of $\Ntilde$ everywhere on the boundary of $\Ntilde$.
That is, for each ${n}\in\Ntilde$ such that ${n}_k=\nbot_k$ (resp. ${n}_k=\ntop_k$) for some $k\in\{1,\ldots,K\}$ (resp. $k\in\{2,\ldots,K\}$), we have $G_k(i,{n},r)\ge0$ (resp. $G_k(i,{n},r)\le0$) for all $i\in \I$.

\emph{a)} We first study the directionality of the vector field on the lower boundaries. Consider a given $r\in\R$.

\emph{a.1)} For each ${n}\in\Ntilde$ such that ${n}_1=\nbot_1$, we have
\begin{align}
G_1(i,{n},r)&\stackrel{\footnotesize\eqref{Eq_G1}}=r_1-f_1(i,\nbot_1,{n}_2,r_2)/\beta_1\nonumber\\
&\stackrel{\footnotesize\eqref{Eq_fk}}=r_1-\min\left\{v\nbot_1,\Sbar_1^i,(R_2({n}_2)-r_2)_+/\beta_2\right\}\nonumber\\
&\ge r_1-v\nbot_1\stackrel{\footnotesize\eqref{Eq_rhobot1}}{\ge}r_1-v\frac{r_1}{v}=0,
\quad\forall i\in\I.
\label{Eq_r1-f1}
\end{align}

\emph{a.2)} For each ${n}\in\Ntilde$ such that ${n}_k=\nbot_k$ for some $k\in\{2,\ldots,K-1\}$, we need to show that $G_k\ge0$.

First, note that
\begin{align}
&f_{k-1}\left(i,{n}_{k-1},\nbot_k,r_k\right)\nonumber\\
&\stackrel{\footnotesize\eqref{Eq_fk}}=\min\left\{\beta_{k-1}v{n}_{k-1},\beta_{k-1}\Sbar_{k-1}^i,R_k(\nbot_k)-r_k\right\}\nonumber\\
&\ge \min\left\{\beta_{k-1}v\nbot_{k-1},\beta_{k-1}\Sbar_{k-1}^i,R_k(\nbot_k)-r_k\right\}\nonumber\\
&\stackrel{\footnotesize\eqref{Eq_Receiving}\eqref{Eq_rhobotk}}{\ge}\min\left\{{v}\nbot_k-r_k,w({n}^\max-\nbot_k)-r_k\right\}
\label{Eq_fk-1}
\end{align}
Since
$
\nbot_k\stackrel{\footnotesize\eqref{Eq_rhobotk}}{\le}
\Sbar_k^\max/v,
$
we can obtain from \eqref{Eq_Sbarmax} that
\begin{align*}
{v}\nbot_k\le w\left({n}^\max-\nbot_k\right).
\end{align*}
Plugging the above into \eqref{Eq_fk-1}, we obtain
\begin{align}
f_{k-1}\left(i,{n}_{k-1},\nbot_k,r_k\right)
\ge{v}\nbot_k-r_k.
\label{Eq_fk-1<=vrho}
\end{align}
Next, note that
\begin{align}
&f_{k}\left(i,\nbot_{k},{n}_{k+1},r_{k+1}\right)/\beta_k\nonumber\\
&\stackrel{\footnotesize\eqref{Eq_fk}}=\min\left\{v\nbot_{k},\Sbar_{k}^i,\frac{w\left({n}^\max-{n}_{k+1}\right)-r_k}{\beta_k}\right\}
\le v\nbot_{k}.
\label{Eq_fk<=vrhok}
\end{align}
Hence, we have
\begin{align}
&G_k(i,{n},k)\nonumber\\
&=f_{k-1}\left(i,{n}_{k-1},\nbot_k,r_k\right)+r_k-f_{k}\left(i,\nbot_{k},{n}_{k+1},r_{k+1}\right)/\beta_k\nonumber\\
&\stackrel{\footnotesize\eqref{Eq_fk-1<=vrho}\eqref{Eq_fk<=vrhok}}\ge{v}\nbot_k-r_k+r_k-v\nbot_k=0,
\quad\forall i\in\I.
\label{Eq_fk-1+rk-fk}
\end{align}

\emph{a.3)} The proof for $k=K$ is analogous.

\emph{b)} Next, we study the directionality of the vector field $G$ on the upper boundaries. Again, consider a given $r\in\R$.

\emph{b.1)} For each ${n}\in\Ntilde$ such that ${n}_K=\ntop_K$, we need to consider the two subcases in \eqref{Eq_rhotopn}:

\underline{If $\beta_{K-1}\Sbar^\max+r_K\le\Sbar_K^\min$}, then we have
\begin{subequations}
	\begin{align}
	&f_{K-1}(i,{n}_{K-1},\ntop_K,r_K)\nonumber\\
	&\stackrel{\footnotesize\eqref{Eq_fk}}=\min\{\beta_{K-1}v{n}_{K-1},\beta_{K-1}\Sbar_{K-1}^i,w({n}^\max-\ntop_{K})-r_K\}\nonumber\\
	&\le\beta_{K-1}F_{K-1}^i\le\beta_{K-1}\Sbar^\max,\label{Eq_fn-1<=betaFmax}\\
	&\frac{f_{K}(i,\ntop_{K})}{\beta_K}
	\stackrel{\footnotesize\eqref{Eq_fn}}=\min\{v\ntop_K,\Sbar_K^i\}
	\stackrel{\footnotesize\eqref{Eq_rhotopn}}{=}\beta_{K-1}\Sbar^\max+r_K.\label{Eq_fn=betaFmax}
	\end{align}%
	\label{Eq_fn-1fn}%
\end{subequations}
	Thus, we have
	\begin{align*}
	&G_K(i,{n},r)
	\stackrel{\footnotesize\eqref{Eq_Gn}}=f_{K-1}(i,{n}_{K-1},\ntop_K,r_K)+r_K-\frac{f_{K}(i_{K},\ntop_{K})}{\beta_K}\\	&\stackrel{\footnotesize\eqref{Eq_fn-1fn}}\le\beta_{K-1}\Sbar^\max+r_K-(\beta_{K-1}\Sbar^\max+r_K)=0,
	\quad\forall i\in\I.
	\end{align*}
	
\underline{Otherwise}, we have
	\begin{subequations}
	\begin{align}
	&f_{K-1}(i,{n}_{K-1},\ntop_K,r_K)
	\stackrel{\footnotesize\eqref{Eq_fk}}\le w({n}^\max-\ntop_K)-r_K\nonumber\\
	&\stackrel{\footnotesize\eqref{Eq_rhotopn}}=\Sbar_K^\min-r_K,\\
	&f_{K}(i_{K},\ntop_{K})/\beta_K
	\stackrel{\footnotesize\eqref{Eq_fn}\eqref{Eq_rhotopn}}=\min\{v({n}^\max-\Sbar_K^\min/w),\Sbar_K^i\}\nonumber\\
	&\stackrel{\footnotesize\eqref{Eq_Sbarmax}}=\Sbar_K^i\ge\Sbar_K^\min.
	\end{align}
	\label{Eq_P1b1}%
	\end{subequations}
	Thus, we have
	\begin{align*}
	&G_K(i,{n},r)\stackrel{\footnotesize\eqref{Eq_Gn}}=f_{K-1}(i,{n}_{K-1},\ntop_K,r_K)+r_K-\frac{f_{K}(i_{K},\ntop_{K})}{\beta_K}\\
	&\stackrel{\footnotesize\eqref{Eq_P1b1}}\le(\Sbar_K^\min-r_K)+r_K-\Sbar_K^\min=0,
	\quad\forall i\in\I.
	\end{align*}

\emph{b.2)} For ${n}\in\Ntilde$ such that ${n}_k=\ntop_k$ for some $k\in\{1,\ldots,K-1\}$, we again need to consider both cases indicated in \eqref{Eq_rhotopk}:

\underline{If $\beta_{k-1}\Sbar^\max+r_k\le\min\{\Sbar_k^\min,\frac{R_{k+1}(\ntop_{k+1})-r_{k+1}}{\beta_{k}}\}$}, then we have
\begin{subequations}
	\begin{align}
	&f_{k-1}(i,{n}_{k-1},\ntop_k,r_k)
	\le\beta_{k-1}F_{k-1}^i
	\le\beta_{k-1}\Sbar^\max,\\
	&f_{k}(i,\ntop_{k},{n}_{k+1},r_{k+1})/\beta_k\nonumber\\
	&\stackrel{\footnotesize\eqref{Eq_fk}}=\min\left\{v\ntop_k,\Sbar_k^i,\frac{w({n}^\max-{n}_{k+1})-r_{k+1}}{\beta_k}\right\}\nonumber\\
	&\ge\min\left\{v\ntop_k,\Sbar_k^i,\frac{w({n}^\max-\ntop_{k+1})-r_{k+1}}{\beta_k}\right\}\nonumber\\
	&\stackrel{\footnotesize\eqref{Eq_rhotopk}}=\min\left\{\beta_{k-1}\Sbar^\max+r_k,\Sbar_k^i,\frac{w({n}^\max-\ntop_{k+1})-r_{k+1}}{\beta_k}\right\}\nonumber\\
	&{=}\beta_{k-1}\Sbar^\max+r_k
	\end{align}
	\label{Eq_b2}%
\end{subequations}
	Thus, we have
	\begin{align*}
	&G_k(i,{n},r)\nonumber\\
	&=f_{k-1}(i_{k-1},{n}_{k-1},\ntop_k,r_k)+r_k-f_{k}(i,\ntop_{k},{n}_{k+1},r_{k+1})/\beta_k\\
	&\stackrel{\footnotesize\eqref{Eq_b2}}\le\beta_{k-1}\Sbar^\max+r_k-(\beta_{k-1}\Sbar^\max+r_k)=0,
	\quad\forall i\in\I.
	\end{align*}
	
\underline{Otherwise}, we have	
\begin{subequations}
	\begin{align}
	&f_{k-1}(i,{n}_{k-1},\ntop_k,r_k)
	\le w({n}^\max-\ntop_k)-r_k\nonumber\\
	&\stackrel{\footnotesize\eqref{Eq_rhotopk}}=\min\left\{\Sbar_k^\min,\frac{w({n}^\max-\ntop_{k+1})-r_{k+1}}{\beta_k}\right\}-r_k,\\
	&f_{k}(i,\ntop_{k},{n}_{k+1},r_{k+1})/\beta_K\nonumber\\
	&=\min\left\{v\ntop_k,\Sbar_k^i,\frac{w({n}^\max-{n}_{k+1})-r_{k+1}}{\beta_k}\right\}\nonumber\\
	&\ge\min\left\{v\ntop_k,\Sbar_k^i,\frac{w({n}^\max-\ntop_{k+1})-r_{k+1}}{\beta_k}\right\}\nonumber\\
	&\stackrel{\footnotesize\eqref{Eq_rhotopk}}{=}\min\left\{\Sbar_k^i,\frac{w({n}^\max-\ntop_{k+1})-r_{k+1}}{\beta_k}\right\}\nonumber\\
	&\ge\min\left\{\Sbar_k^\min,\frac{w({n}^\max-\ntop_{k+1})-r_{k+1}}{\beta_k}\right\}.
	\end{align}
	\label{Eq_b2otherwise}
\end{subequations}
	Thus, we have
	\begin{align*}
	&G_k(i,{n},r)=f_{k-1}(i,{n}_{k-1},\ntop_k,r_k)+r_k\nonumber\\
	&\quad-f_{k}(i_{k},\ntop_{k},{n}_{k+1},r_{k+1})/\beta_k\stackrel{\footnotesize\eqref{Eq_b2otherwise}}\le0,
	\quad\forall i\in\I.
	\end{align*}

Combining cases a) and b), we obtain that $\Ntilde$ is invariant.

\subsubsection{Attracting} To show that the set $\Ntilde$ is attracting, we define
\begin{align*}
&\underline\N_k=\{{n}\in\N:{n}_k\ge\nbot_k\},\; k=1,2,\ldots,n,\\
&\overline\N_k=\{{n}\in\N:{n}_k\le\ntop_k\},\; k=2,3,\ldots,n.
\end{align*}
Thus, we have $\Ntilde=(\cap_{k=1}^K\underline\N_k)\cap(\cap_{k=2}^K\overline\N_k)$.
Consider a given $r\in\R$.

\emph{a)} First, we show by induction that the set $\cap_{k=1}^K\underline\N_k$ is attracting.

\emph{a.1)} For any $\epsilon>0$, consider $\B(\underline\N_1,\epsilon)$, i.e. the neighborhood of $\underline\N_1$ such that $\min_{\varrho\in\underline\N_1}\|{n}-\varrho\|_2\le\epsilon$.
Without loss of generality, we consider $0<\epsilon<\min_{k:\nbot_k>0}\nbot_k$.
If $\nbot_1>0$, for any ${n}\in\B^{\mathrm c}(\underline\N_1,\epsilon)$ (complement of $\B(\underline\N_1,\epsilon)$), we have ${n}_1<\nbot_1-\epsilon$. Then, we obtain from \eqref{Eq_Sending} and \eqref{Eq_f(rho)} that
\begin{align*}
&G_1(i,{n},r)=r_1-f_1(i,{n}_1,{n}_2,r_2)/\beta_1
\ge r_1-v{n}_1\nonumber\\
&
\quad\ge r_1-v(\nbot_1-\epsilon)
\ge v\epsilon>0,
	\quad\forall i\in\I.
\end{align*}
Therefore, for any initial condition $(i,{n})\in \I\times\B^{\mathrm c}(\underline\N_1,\epsilon)$, there exists $T=\nbot_1/(v\epsilon)$ such that ${N}(t)\in\B(\underline\N_1,\epsilon)$ for all $t\ge T$.
Hence, the set $\underline\N_1$ is attracting in the sense of \eqref{Eq_Attracting}.

If $\nbot_1=0$, the proof is trivial.

\emph{a.2)} Now, suppose that the set $(\cap_{h=1}^{k}\underline\N_h)$ is attracting.
If $\nbot_{k+1}=0$, for any $\epsilon>0$, consider the neighborhood $\B(\cap_{h=1}^{k+1}\underline\N_h,\epsilon)$.
For each ${n}\in(\cap_{h=1}^{k}\underline\N_h)\cap\B^{\mathrm c}(\cap_{h=1}^{k+1}\underline\N_h,\epsilon)$, we have ${n}_1\ge\nbot_1,\ldots,{n}_k\ge\nbot_k,{n}_{k+1}<\nbot_{k+1}-\epsilon$.
Then, we obtain from \eqref{Eq_Sending} and \eqref{Eq_f(rho)} that
\begin{align*}
&G_{k+1}(i,{n},r)=f_k(i,{n}_k,{n}_{k+1},r_{k+1})+r_{k+1}\nonumber\\
&\quad\quad-f_{k+1}(i,{n}_{k+1},{n}_{k+2},r_{k+2})/\beta_{k+1}\\
&\quad\stackrel{\footnotesize\eqref{Eq_fk}}\ge f_k\left(i,\nbot_k,{n}_{k+1},r_{k+1}\right)+r_{k+1}\nonumber\\
&\quad\quad-f_{k+1}(i,{n}_{k+1},{n}_{k+2},r_{k+2})/\beta_{k+1}\\
&\quad\stackrel{\footnotesize\eqref{Eq_rhobotk}}\ge(\nbot_{k+1}v-r_{k+1})+r_{k+1}-v\left(\nbot_{k+1}-\epsilon\right)\\
&\quad\ge v\epsilon>0,
	\quad\forall i\in\I.
\end{align*}
Therefore, for any initial condition $(i,{n})\in \I\times(\cap_{h=1}^{k}\underline\N_h)\cap\B^{\mathrm c}(\cap_{h=1}^{k+1}\underline\N_h,\epsilon)$, there exists $T=\nbot_{k+1}/(v\epsilon)$ such that ${N}(t)\in\B(\cap_{h=1}^{k+1}\underline\N_h,\epsilon)$ for all $t\ge T$.
Recall that, by the inductive hypothesis, $\cap_{h=1}^{k}\underline\N_h$ is (globally) attracting.
Hence, the set $\cap_{h=1}^{k+1}\underline\N_h$ is attracting.

If $\nbot_{k+1}=0$, the proof is trivial.

In conclusion, $\cap_{k=1}^K\underline\N_k$ is attracting.

\emph{b)} The proof for $\cap_{k=2}^K\overline\N_k$ being attracting is analogous.

\subsection{Proof of Theorem~\ref{Thm_Stable1}}

Suppose that the SS-CTM with a given inflow vector $r$ is stable in the sense of \eqref{Eq_Stable}.
To establish the necessary condition, we can limit our attention to a particular initial condition in the invariant set, i.e. $N(0)={n}\in\Ntilde$.
The proof consists of two steps. In Step 1), we show that the average flow is equal to the nominal flow. In Step 2), we show that the average flow is less than or equal to the average spillback-adjusted capacity.

\emph{1)} Integrating \eqref{Eq_rhodot}, we obtain that, for $t\ge0$,
\begin{align*}
&{N}_k(t)=\int_{\tau=0}^t\Big(f_{k-1}(\tau)+r_k-f_k(\tau)/\beta_k\Big)d\tau+{n}_k,\\
&\hspace{2in}k=1,2\ldots,n.
\end{align*}
Since $\lim_{t\to\infty}{n}_k/t=0$ for $k=1,2,\ldots,n$, we can write
\begin{align*}
0&=\lim_{t\to\infty}\frac{1}{t}\Bigg(\int_{\tau=0}^t\Big(f_{k-1}(\tau)+r_k-f_k(\tau)/\beta_k\Big)d\tau\nonumber\\
&\hspace{.8in}+{n}_k-{N}_k(t)\Bigg)\nonumber\\
&=\lim_{t\to\infty}\frac{1}{t}\Bigg(\int_{\tau=0}^t\Big(f_{k-1}(\tau)+r_k-f_k(\tau)/\beta_k\Big)d\tau-{N}_k(t)\Bigg).
\end{align*}
Since the MGF of $|N(t)|$ is bounded on average, we have $\Pr\{\lim_{t\to\infty}{N}_k(t)=\infty\}=0$ for $k=1,2,\ldots,n$.
Thus, we have
\begin{align}
\lim_{t\to\infty}\frac{1}{t}\int_{\tau=0}^t\Big(f_{k-1}(\tau)+r_k-\frac{f_k(\tau)}{\beta_k}\Big)d\tau=0,\quad a.s.
\label{Eq_limphi}
\end{align}

For $k=1$, since $f_0=0$ by definition, we have
\begin{align}
\lim_{t\to\infty}\frac{1}{t}\int_{\tau=0}^t\Big(r_1-f_1(\tau)/\beta_1\Big)d\tau=0,\quad a.s.
\label{Eq_rho1}
\end{align}
which implies that
\begin{align*}
\lim_{t\to\infty}\frac1t\int_{\tau=0}^tf_1(\tau)d\tau=\beta_1r_1.\quad a.s.
\end{align*}
To proceed by induction, we assume that, for some $k\ge1$, we have
\begin{align}
\lim_{t\to\infty}\frac1t\int_{\tau=0}^tf_k(\tau)d\tau=\beta_k\sum_{h=1}^k\beta_h^kr_h,\quad a.s.
\label{Eq_limfk}
\end{align}
Then, we can obtain from \eqref{Eq_limphi} that
\begin{align*}
&\lim_{t\to\infty}\frac1t\int_{\tau=0}^tf_{k+1}(\tau)d\tau\\
&=\beta_{k+1}\left(\lim_{t\to\infty}\frac1t\int_{\tau=0}^tf_k(\tau)d\tau+r_{k+1}\right)\\
&\stackrel{\footnotesize\eqref{Eq_limfk}}=\beta_{k+1}\left(\beta_k\sum_{h=1}^k\beta_h^kr_h+r_{k+1}\right)
=\beta_{k+1}\sum_{h=1}^{k+1}\beta_h^{k+1}r_h,\; a.s.
\end{align*}

Hence, we conclude that
\begin{align}
\lim_{t\to\infty}\frac1t\int_{\tau=0}^tf_k(\tau)d\tau=\beta_k\sum_{h=1}^k\beta_h^kr_h,\; a.s.\; k=1,2,\ldots,n,
\label{Eq_lim1}
\end{align}
which means that the average flow is equal to the nominal flow.

\emph{2)} For every $i\in\I$, let $T_i(t)$ be the amount of time that the SS-CTM is in mode $i$ up to time $t$, i.e.
\begin{align}
T_i(t)=\int_{\tau=0}^t\mathbb1_{\{I(\tau)=i\}}d\tau.
\label{Eq_Ti}
\end{align}
Then, recalling Assumption~\ref{Asm_Ergodic} and using \cite[Theorem 7.2.6]{gallager13}, we obtain
\begin{align}
\lim_{t\to\infty}\frac{T_i(t)}t=\sfp_i,\quad a.s.
\label{Eq_T/t}
\end{align}
In addition, since $\Ntilde$ is invariant, we know from \eqref{Eq_Invariant0} that ${n}(\tau)\in\Ntilde$ for all $\tau\ge0$. Then, we have, for $k=1,2,\ldots,n$,
\begin{align}
&\lim_{t\to\infty}\frac{1}{t}\int_{\tau=0}^tf_k(\tau)d\tau
\stackrel{\footnotesize\eqref{Eq_fk<=min}}{\le}\lim_{t\to\infty}\frac{1}{t}\int_{\tau=0}^t\beta_k\Stilde_k^{I(\tau)}(\nbot,r)d\tau\nonumber\\
&\stackrel{\footnotesize\eqref{Eq_Ti}}=\beta_k\lim_{t\to\infty}\sum_{i\in\I}\frac{T_i(t)\Stilde_k^i(\nbot,r)}{t}
\stackrel{\footnotesize\eqref{Eq_T/t}}=\beta_k\sum_{i\in\I}\sfp_i\Stilde_k^i(\nbot,r).\quad a.s.
\label{Eq_lim2}
\end{align}

Combining \eqref{Eq_lim1} and \eqref{Eq_lim2}, we obtain \eqref{Eq_r<=Stilde}.

\subsection{Proof of Theorem~\ref{Thm_Stable2}}
\label{Sub_Stable2}

Recall from Section~\ref{Sub_Results} that, to obtain stability, we need to show \eqref{Eq_LV2}.
Let us proceed with the following expression:
\begin{align}
&\max_{{n}\in\Ntilde}\Ell V(i,{n})+cV(i,{n})\nonumber\\
&\stackrel{\footnotesize\eqref{Eq_Lg}\eqref{Eq_V}}=\max_{{n}\in\Ntilde}\Bigg(a_ib\sum_{k=1}^K\Gamma_k(f_{k-1}+r_k-f_k/\beta)\nonumber\\
&\hspace{.9in}+\sum_{j=1}^m\lambda_{ij}(a_j-a_i)+a_ic\Bigg)\frac{V}{a_i}\nonumber\\
&\stackrel{\footnotesize\eqref{Eq_Wk}\eqref{Eq_scrR}}=\Bigg(a_ib\Big(\scrR(r)-\min_{{n}\in\Ntilde}\gamma^Tf(i,{n},r)\Big)\nonumber\\
&\hspace{.6in}+\sum_{j=1}^m\lambda_{ij}(a_j-a_i)+a_ic\Bigg)\frac{V}{a_i},\;
\forall i\in \I,
\label{Eq_LV3}
\end{align}
where $c$ and $d$ are given by \eqref{Eq_cd}.
The key to evaluate \eqref{Eq_LV3} is to compute $\min_{{n}\in\Ntilde}\gamma^Tf(i,{n},r)$.
To do this, we define
\begin{align*}
&\Ntilde_1:=\left[{n}^{\mathrm{crit}},\infty\right)\times\prod_{k=2}^K\left[\nbot_k,\ntop_k\right],\\
&\Ntilde_2:=\left[\nbot_1,{n}^{\mathrm{crit}}\right]\times\prod_{k=2}^K\left[\nbot_k,\ntop_k\right]
\end{align*}
and consider two cases:

\emph{1) ${n}\in\Ntilde_1$.}
For each $r\in\R$ and each $i\in \I$, consider
\begin{align}
\min_{{n}\in\Ntilde_1}\quad \gamma^Tf(i,{n},r).\label{Eq_minF}
\end{align}
We claim that an optimal solution of \eqref{Eq_minF} lies on one of the vertices of $\Ntilde_1$, i.e. for all $i\in \I$,
\begin{align}
\min_{{n}\in\Ntilde_1}\gamma^Tf(i,{n},r)
=\min_{{n}\in\Theta}\gamma^Tf(i,{n},r)
\stackrel{\small\eqref{Eq_Fi}}=\scrF_i(\nbot,\ntop,r),
\label{Eq_Rhotilde=Theta}
\end{align}
where $\Theta$ is defined in \eqref{Eq_Theta}. We will prove this claim at the end of this subsection.

Then, for each $i\in \I$, we have
\begin{align*}
&\max_{{n}\in\Ntilde_1}\Ell V(i,{n})+cV(i,{n})\\
&\stackrel{\footnotesize\eqref{Eq_LV3}\eqref{Eq_Rhotilde=Theta}}\le\Bigg(a_ib\left(\scrR(r)-\scrF_i(\nbot,\ntop,r)\right)+\sum_{j=1}^m\lambda_{ij}(a_j-a_i)+a_ic\Bigg)\\
&\hspace{.4in}\times\exp\left(\Gamma^T{n}\right)\\
&\stackrel{\footnotesize\eqref{Eq_BMI}}\le(-1+a_ic)\exp\left(\Gamma^T{n}\right)\stackrel{\footnotesize\eqref{Eq_c}}\le0\stackrel{\footnotesize\eqref{Eq_d}}\le d.
\end{align*}

\emph{2) ${n}\in\Ntilde_2$.}
This case is straightforward, since $\Ntilde_2$ is bounded.
We claim that, for all $i\in \I$,
\begin{align}
\min_{{n}\in\Ntilde_2}\gamma^Tf(i,{n},r)
=\min_{{n}\in\hat\Theta}\gamma^Tf(i,{n},r)
\stackrel{\small\eqref{Eq_Fihat}}=\hat\scrF_i(\nbot,\ntop,r),
\label{Eq_gammaf>scrFhat}
\end{align}
where $\hat\Theta$ is as defined in \eqref{Eq_Thetahat};
again, we will prove this claim at the end of this subsection.
Then, we have
\begin{align*}
&\max_{{n}\in\Ntilde_2}\Ell V+cV\\
&\le\Bigg|a_ib\left(\scrR(r)-\hat\scrF_i(\nbot,\ntop,r)\right)+\sum_{j=1}^m\lambda_{ij}(a_j-a_i)+a_ic\Bigg|\nonumber\\
&\quad\times\frac{\max_{{n}\in\Ntilde_2}V}{a_i}
\stackrel{\footnotesize\eqref{Eq_d}}\le d,\; \forall i\in \I.
\end{align*}

Since $\Ntilde=\Ntilde_1\cup\Ntilde_2$, combining cases a) and b), we obtain \eqref{Eq_LV2} and thus the drift condition \eqref{Eq_Drift0}. Using \cite[Theorem 4.3]{meyn93}, we obtain \eqref{Eq_Stable2} from the drift condition.

The rest of this subsection is devoted to the derivation of \eqref{Eq_Rhotilde=Theta} and \eqref{Eq_gammaf>scrFhat}:

To show \eqref{Eq_Rhotilde=Theta}, it suffices to argue that, for a given $i\in \I$ and a given $r\in\R$, and for each $k\in\{1,\ldots,K\}$, $\gamma^Tf(i,{n},r)$ is concave in ${n}_k$.
For each $i\in \I$, let us consider the following quantity:
\begin{align*}
H_k^i({n})=\left\{\begin{array}{ll}
\gamma_1f_1(i,{n}_1,{n}_2,r_2), & k=1,\\
\gamma_{k-1}f_{k-1}(i,{n}_{k-1},{n}_k,r_k) & \\
\quad+\gamma_{k}f_{k}(i,{n}_{k},{n}_{k+1},r_{k+1}), & 2\le k\le K-1,\\
\gamma_{K-1}f_{K-1}(i,{n}_{K-1},{n}_K,r_K)&\\
\quad+\gamma_{K}f_{K}(i,{n}_{K}), &k=K.
\end{array}\right.
\end{align*}
Then, for each $k\in\{1,\ldots,K\}$, we have
\begin{align*}
\gamma^Tf(i,{n},r)=H_k^i({n})+M_k^i({n}_1,\ldots,{n}_{k-1},{n}_{k+1},\ldots,{n}_K),
\end{align*}
where $M_k^i$ is a term independent of ${n}_k$.
Hence, to show that $\gamma^Tf(i,{n},r)$ is concave in ${n}_k$, we only need to show that $H_k({n})$ is concave in ${n}_k$.

\begin{figure}[hbt]
\centering
\includegraphics[width=0.3\textwidth]{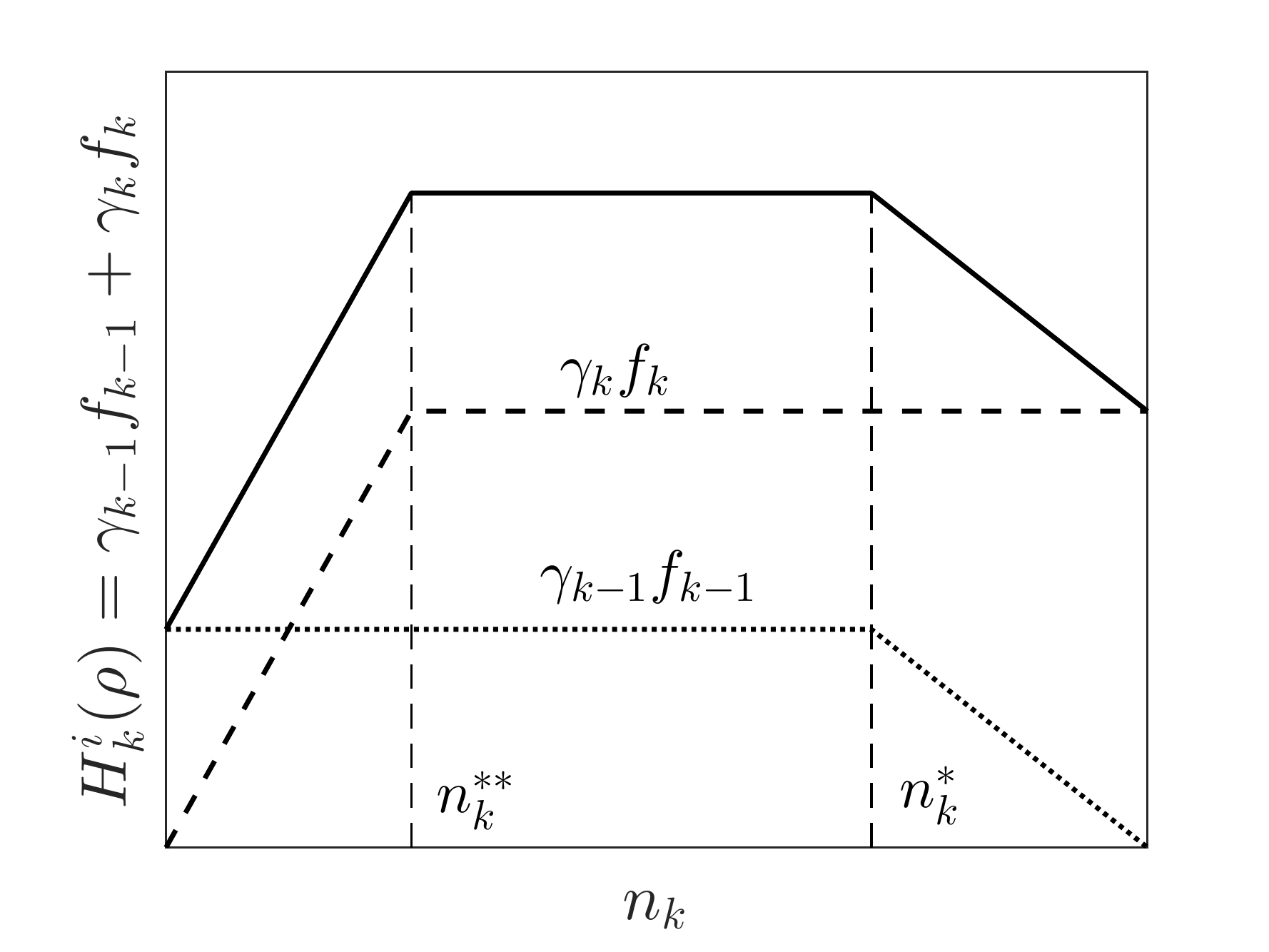}
\caption{The function $H_k({n})$ is concave in ${n}_k$.}
\label{Eq_wfvrhok}
\end{figure}

We need to consider the following subcases of $k$:

\emph{a.1)}:
For each $k\in\{2,\ldots,K-1\}$, we have
\begin{align*}
&H_k({n})\\
&=\gamma_{k-1}\min\left\{\beta_{k-1}v{n}_{k-1},\beta_{k-1}\Sbar_{k-1}^i,w({n}^\max-{n}_k)-r_k\right\}\\
&\quad+\gamma_{k}\min\left\{\beta_{k}v{n}_{k},\beta_{k}\Sbar_{k}^i,w({n}^\max-{n}_{k+1})-r_{k+1}\right\}.
\end{align*}
In the above, the first term on the right side (corresponding to $\gamma_{k-1}f_{k-1}$) is non-increasing in ${n}_k$, while the second term (corresponding to $\gamma_kf_k$) is non-decreasing in ${n}_k$; both terms are piecewise affine in ${n}_k$ with exactly one intersecting point (see Figure~\ref{Eq_wfvrhok}).

Note that the intersecting points ${n}_k^*$ and ${n}_k^{**}$ of the piecewise-linear functions $\gamma_{k-1}f_{k-1}$ and $\gamma_kf_k$, respectively, satisfy the following:
\begin{align*}
{n}_k^*&={n}^\max-\frac1w(\min\{\beta_{k-1}v{n}_{k-1},\beta_{k-1}F_{k-1}^i\}-r_k)\\
&\ge{n}^\max-\frac{\beta_{k-1}F_{k-1}^i}{w}
\ge{n}^\max-\frac{F^\max}{w}={n}^{\mathrm{crit}},\\
{n}_k^{**}&=\frac1{\beta_kv}\min\{\beta_{k}F_{k}^i,w({n}^\max-{n}_{k+1})-r_{k+1}\}\\
&\le\frac{F_{k}^i}{v}\le\frac{F^\max}{v}
={n}^{\mathrm{crit}},
\end{align*}
Hence, we have ${n}_k^{*}\ge {n}_k^{**}$.
Thus, we conclude that $H_k({n})$ is concave in ${n}_k$; see Figure~\ref{Eq_wfvrhok}.


\emph{a.2)}:
For $k=1$ and $k=K$, the expression of $H_k$ is simpler and thus the derivation is straightforward.

The proof of \eqref{Eq_gammaf>scrFhat} is analogous.

\subsection{Proof of Proposition~\ref{Prp_Rhotilde}}

We show the two conclusions in the statement separately:

\emph{(i)}: We can observe from \eqref{Eq_Ftilde} that $\Ftilde_k^i(\nbot,r)$ is non-increasing in $\nbot$. Hence, for each $\nbot'\le\nbot$, we have
\begin{align*}
\Ftilde_k^i(r,\nbot')\ge\Ftilde_k^i(\nbot,r),
\quad k=1,2,\ldots,n,\;i\in \I.
\end{align*}
The conclusion follows from the above.

\emph{(ii)}: Consider $\Ntilde\subseteq\Ntilde'$. Define
\begin{align*}
&\Ntilde_1':=\left[{n}^{\mathrm{crit}},\infty\right)\times\prod_{k=2}^K\left[\nbot_k',\ntop_k'\right],\\
&\Ntilde_2':=\left[\nbot_1',{n}^{\mathrm{crit}}\right]\times\prod_{k=2}^K\left[\nbot_k',\ntop_k'\right].
\end{align*}
Clearly, $\Ntilde_1\subseteq\Ntilde_1'$ and $\Ntilde_2\subseteq\Ntilde_2'$.
Then, from \eqref{Eq_Rhotilde=Theta} and \eqref{Eq_gammaf>scrFhat}, we can obtain that
\begin{align*}
&\scrF_i(\nbot,\ntop,r)=\min_{{n}\in\Ntilde_1}\gamma^Tf(i,{n},r)\ge
\min_{{n}'\in\Ntilde_1'}\gamma^Tf(i,{n}',r)\\
&=\scrF_i(r,\nbot',\ntop'),\quad\forall i\in \I,\\
&\scrF_i(\nbot,\ntop,r)=\min_{{n}\in\Ntilde_2}\gamma^Tf(i,{n},r)\ge
\min_{{n}'\in\Ntilde_2'}\gamma^Tf(i,{n}',r)\\
&=\scrF_i(r,\nbot',\ntop'),\quad\forall i\in \I,
\end{align*}
The conclusion follows from the above and the fact that $a_1,a_2,\ldots,a_m$ and $b$ in \eqref{Eq_BMI} are positive.
%
\section*{Acknowledgment}
This work was supported by NSF Grant No. 1239054 CPS Frontiers (FORCES), MIT CEE Graduate Fellowship, and NSF CAREER Award CNS 1453126. We also sincerely appreciate the valuable inputs from the associate editor and the anonymous reviewers.

\bibliographystyle{plain}
\bibliography{Bibliography}

\end{document}